\documentclass{amsart}

\usepackage{amssymb,amsthm,amsfonts,latexsym,mathtools}
\usepackage[utf8]{inputenc}
\usepackage{xcolor}
\usepackage{mathtools}
\usepackage{hyperref}
\usepackage{mathdots}
\usepackage{array}
\usepackage{enumitem}

\newcommand{\C}{\mathbb{C}}
\newcommand{\F}{\mathbb{F}}

\newcommand{\calC}{\mathcal{C}}
\newcommand{\calG}{\mathcal{G}}
\newcommand{\calL}{\mathcal{L}}

\newcommand{\Aut}{\mathrm{Aut}}

\newcommand{\Sym}{\mathbb{S}}
\newcommand{\Alt}{\mathbb{A}}
\newcommand{\SL}{\mathrm{SL}}
\newcommand{\PSL}{\mathrm{PSL}}
\newcommand{\PGL}{\mathrm{PGL}}
\newcommand{\PSU}{\mathrm{PSU}}
\newcommand{\PGU}{\mathrm{PGU}}
\newcommand{\SU}{\mathrm{SU}}
\newcommand{\GU}{\mathrm{GU}}
\newcommand{\Sz}{\mathrm{Sz}}
\newcommand{\Ree}{\mathrm{Ree}}

\newcommand{\GL}{\mathrm{GL}}

\newcommand{\Syl}{\mathrm{Syl}}
\newcommand{\Irr}{\mathrm{Irr}}
\newcommand{\normal}{\trianglelefteq}

\newcommand{\gen}[1]{\langle #1\rangle}
\newcommand{\Tr}{\mathrm{Tr}}
\newcommand{\ad}{\mathrm{ad}}

\newcommand{\commGraph}[1]{\calG(#1)}
\newcommand{\killingGraph}[2]{\calG_K(#1, #2)}

\numberwithin{equation}{section}

\newtheorem{theorem}{Theorem}[section]
\newtheorem{lemma}[theorem]{Lemma}
\newtheorem{proposition}[theorem]{Proposition}
\newtheorem{corollary}[theorem]{Corollary}
\newtheorem{remark}[theorem]{Remark}
\newtheorem{question}[theorem]{Question}
\newtheorem{conjecture}[theorem]{Conjecture}

\theoremstyle{definition}
\newtheorem{definition}[theorem]{Definition}
\newtheorem{example}[theorem]{Example}

\title[On reducible Killing forms for groups of Lie type]{On reducible Killing forms for groups of Lie type}

\author{Kevin I. Piterman}
\email{kevin.piterman@vub.be}

\author{Charlotte Roelants}
\email{charlotte.roelants@vub.be }

\address{Department of Mathematics and Data
Science, Vrije Universiteit Brussel, Pleinlaan 2, 1050 Brussel, Belgium}


\begin{document}

\begin{abstract}
Killing forms on finite groups arise as examples of braided Killing forms on braided Lie algebras.
For a finite group $G$ and a $G$-stable subset $\calC$, the Killing form associated with $\mathbb{C}[\calC]$ is given by $K_\calC(a,b) = |C_G(ab) \cap \calC|$ for $a,b\in \calC$.  
Motivated by Cartan's criterion for semisimplicity of Lie algebras, and previous work of L{\'o}pez Pe{\~n}a, Majid, and Rietsch, we study the non-degeneracy and irreducibility of $K_{\calC}$ when $\calC$ is a conjugacy class of involutions or unipotent elements in a finite simple group of Lie type and Lie rank one.
Our approach suggests interesting connections with character theory, related counting formulas, and the study of commuting graphs.
\end{abstract}

\maketitle

\section{Introduction}

Killing forms were originally introduced in the context of Lie algebras: for $\mathfrak{g}$ a Lie algebra over the field $k$ with Lie bracket $[\, \cdot \, , \, \cdot \, ]$, the \emph{Killing form} associated with $\mathfrak{g}$ is the bilinear form 
\[
K: \mathfrak{g} \times \mathfrak{g} \longrightarrow k, \quad 
K(x,y) = \Tr(\ad(x) \circ \ad(y)),
\]
where $\ad(x): \mathfrak{g} \longrightarrow \mathfrak{g}$, given by $z \mapsto [x,z]$ denotes the adjoint endomorphism. 
This form is said to be \emph{non-degenerate} if it is invertible, and it is \textit{degenerate} otherwise.
By Cartan's criterion, a finite-dimensional Lie algebra over a field of characteristic zero is semisimple if and only if its Killing form is non-degenerate (see \cite{Cartan}).

Generalising this construction, in \cite{Majid_braided}, Majid introduced the notion of braided Lie algebras $(\mathcal{L}, \Delta, \epsilon, [ \, \cdot \, , \cdot \, ])$, and defined its corresponding braided Killing form as the braided trace of the map $[\, \cdot \, , \, [\, \cdot \, , \, \cdot \, ] \, ]$.
Similarly, the Killing form is said to be \textit{non-degenerate} if it is invertible.
This motivates the study of an analogue of Cartan's criterion in the broader context of braided Lie algebras.

As proposed in \cite{MR3682632}, when $\calL$ is the $\C$-vector space spanned by certain elements of a finite group, there is a natural braided Lie algebra structure associated with $\calL$ that yields a Killing form whose values are closely related to character computations.
Concretely, let $G$ be a finite group and $\calC \subseteq G \setminus\{1\}$ a $G$-stable subset. 
Take $\calL = \C \calC$, which can be viewed as an object in the braided category of $\C$-vector spaces with trivial braiding operator and unit object $\C$. 
The maps $\Delta$, $\epsilon$ and $[\, \cdot \, , \, \cdot \, ]$ are defined by setting $\Delta(a) = a \otimes a$, $\epsilon(a) = 1$ for $a \in \calC$, and $[a,b] = aba^{-1}$ for $a,b \in \calC$. 
One can verify that, with these choices, $(\mathcal{L}, \Delta, \epsilon, [ \, \cdot \, , \cdot \, ])$ becomes a braided Lie algebra (see \cite{MR3682632,Majid_braided} for more details).
The Killing form associated with $\calC$ is that of $(\mathcal{L}, \Delta, \epsilon, [ \, \cdot \, , \cdot \, ])$, and it is given by
\[
K_\calC(a,b) = \Tr_\calL([a, [b, \, \cdot \, ]]) = |\{ x \in \calC : abxb^{-1}a^{-1} = x \}| = |C_G(ab) \cap \calC|.
\]
The trace in the second expression above yields an interesting connection to character theory, as the Killing form can also be written as $K_\calC(a,b) = \chi(ab)$, with $\chi$ the character of the conjugation representation
\[
\rho: G \longrightarrow \GL(\C \calC), \quad \rho_g(x) = gxg^{-1}.
\]
Moreover, if $\calC = G\setminus \{1\}$, by \cite[Theorem 4.2]{MR3682632}, $K_{\calC}$ is non-degenerate if $G$ has the \textit{Roth property}.
The latter means that every irreducible representation is contained in the conjugation representation of $G$ on itself.
Indeed, \cite{roth_simplegroups} proves that the only finite simple groups failing the Roth property are the unitary groups $\PSU_n(q)$ with $n$ prime to $2(q+1)$, and, in such cases, exactly one irreducible is missing.
Also, \cite{MR3682632} shows that for $G = \PSU_3(3)$, the Killing form associated with $\calC = G\setminus \{1\}$ is degenerate.

Supported by computer calculations for non-abelian simple groups up to order 75000, in \cite{MR3682632} the authors raise the following conjecture:

\begin{conjecture}
\label{conj:class-non-degen}
For every finite non-abelian simple group, the Killing form $K_{\calC}$ on a real conjugacy class $\calC$ is non-degenerate.
\end{conjecture}

So far, there is only little evidence on Conjecture \ref{conj:class-non-degen}.
In \cite{MR3682632}, it is proved that for the conjugacy class of transpositions in the symmetric group, its corresponding Killing form is non-degenerate.
However, to our knowledge, there is no infinite family of non-abelian simple groups satisfying Conjecture \ref{conj:class-non-degen}, not even for specific families of conjugacy classes.

In view of the connection with the Roth property, it would be interesting to obtain a ``$\calC$-local" version of the Roth property for more general $G$-stable subsets $\calC$ of $G\setminus \{1\}$, which implies the non-degeneracy of its corresponding Killing form and hence gives sufficient conditions in terms of irreducible representations to establish Conjecture \ref{conj:class-non-degen}.

On the other hand, following \cite{MR3682632}, when the underlying matrix of $K_{\calC}$ admits no non-trivial diagonal-block decomposition in terms of some ordering of $\calC$,
then $K_{\calC}$ is termed \textit{irreducible}.
Otherwise, $K_{\calC}$ is said to be \textit{reducible}.
Indeed, the latter can be formulated equivalently in terms of the connectivity of a given graph.
For $\calC \subseteq G \setminus  Z(G)$ a $G$-stable subset, let $\killingGraph{G}{\calC}$ be the graph whose vertices are the elements of $\calC$, and with edges between vertices $a,b\in \calC$ if $K_{\calC}(a,b)\neq 0$, that is, if there is an element $c\in \calC$ that commutes with the product $ab$.
Then, $K_{\calC}$ is irreducible if and only if $\killingGraph{G}{\calC}$ is connected  (see Section \ref{sec:overall} for more details).
Observe that the commuting graph with vertices $\calC$ and edges joining commuting vertices, which we denote by $\commGraph{\calC}$, is a subgraph of $\killingGraph{G}{\calC}$ containing all the vertices. In particular, if $\commGraph{\calC}$ is connected, then so is $\killingGraph{G}{\calC}$ and hence $K_{\calC}$ is irreducible.
The study of the connectivity and other combinatorial properties of commuting graphs has attracted the attention of numerous group theorists (see \cite{commGraph1, commGraph2, MR1994533, commGraph3, commGraph4}).

In \cite{MR3682632}, the authors computationally verify that, for non-abelian simple groups $G$ up to order 75000, if $\killingGraph{G}{\calC}$ is disconnected for $\calC$ a conjugacy class of $G$, then $G$ has a strongly embedded subgroup and $\calC$ is the (unique) conjugacy class of involutions.
Groups with a strongly embedded subgroup are classified by the Bender-Suzuki theorem \cite{Bender}, and, for simple groups, these are exactly the groups of Lie type and Lie rank one in characteristic two (these are usually called ``Bender groups").
Equivalently, these are the simple groups for which the commuting graph with vertices all the involutions is disconnected (see \cite{Aschbacher1}).

In Theorem \ref{thm:involutions_Lietype}, we show that for simple groups with a strongly embedded subgroup and $\calC$ the class of involutions, $K_{\calC}$ is reducible.
Moreover, in Theorem \ref{thm:summarySymn_Altn}, we show that for alternating groups $G = \Alt_n$ with $n\equiv 1\pmod 4$, and $\calC$ the conjugacy class of involutions with exactly one fixed point, the Killing form $K_{\calC}$ is reducible, providing a new family of examples.
Nevertheless, we have found no example so far of a reducible Killing form for a conjugacy class of non-involutive elements in a simple group.
Motivated by these computations, we raise the following question:

\begin{question} \label{ques:reduc}
If $G$ is a finite non-abelian simple group and $\calC\subseteq G\setminus \{1\}$ is a $G$-conjugacy class such that $K_{\calC}$ is reducible, does $\calC$ consist of involutions?
\end{question}

We note that, if $G$ is not a simple group, $K_{\calC}$ can be reducible while $\calC$ is not a conjugacy class of involutions.
See Example \ref{ex:reducible_order3}.

From a group-theoretic perspective, it is also interesting to study the structure and the role of the stabiliser of a connected component of the graph $\killingGraph{G}{\calC}$ when this is disconnected.
For instance, for Bender groups $G$ and $\calC \subseteq G$ the conjugacy class of involutions, we will see that $\killingGraph{G}{\calC}$ equals the commuting graph on $\calC$, and so the stabiliser of a component is a strongly embedded subgroup, which is a Borel subgroup.
For alternating groups $G = \Alt_n$ of degree $n\equiv 1\pmod 4$ with $\calC$ the conjugacy class of involutions with one fixed point, $\killingGraph{G}{\calC}$ also equals the commuting graph on $\calC$, which is totally disconnected and hence the stabiliser of a component, i.e., a vertex $x\in \calC$, is just the centraliser $C_G(x)$.
Remarkably, in these cases where $K_{\calC}$ is reducible, it is also non-degenerate, providing an infinite family of examples satisfying Conjecture \ref{conj:class-non-degen} at a given type of conjugacy class.
See Theorem \ref{thm:involutions_Lietype} below.

A notable family to study next is that of the conjugacy classes of $p$-elements in simple groups with a strongly $p$-embedded subgroup, for any prime $p$.
In these cases, the corresponding commuting graphs are always disconnected, and, for $p$-rank at least two, these groups are classified (see Theorem \ref{thm:stronglypembedded}).
Recall that the $p$-rank of a finite group $G$ is the largest rank of an elementary abelian $p$-subgroup of $G$.
In contrast to the case $p=2$ with involutions, for odd $p$ we will see that the associated Killing form on a conjugacy class of non-trivial $p$-elements is always irreducible.

Below, we summarise the main results of this article.

\begin{theorem} \label{thm:involutions_Lietype}
Let $G$ be a finite simple group of Lie type and Lie rank one, i.e., one of the groups
\[ \PSL_2(q) \, (q\geq 4), \ \PSU_3(q)\, (q\neq 2), \ \Sz(q) \, (q\neq 2), \ \Ree(q) \, (q\neq 3).\]
Let $\calC$ be the conjugacy class of involutions in $G$.
Then $K_\calC$ is reducible if and only if $q$ is even. In this case, $\killingGraph{G}{\calC}$ equals the commuting graph on $\calC$, so the number of connected components coincides with the number of Sylow $2$-subgroups (with stabiliser given by a Borel subgroup), and $K_{\calC}$ is non-degenerate.
\end{theorem}

This theorem then gives us an infinite family of groups for which Conjecture \ref{conj:class-non-degen} holds in the case $\calC$ is a conjugacy class of involutions.

Moreover, for simple groups with a strongly $p$-embedded subgroup and $p$-rank at least two and $\calC$ a class of non-trivial $p$-elements, we show that $K_{\calC}$ is irreducible except in the cases of Theorem \ref{thm:involutions_Lietype}.
This leads to the following result:

\begin{theorem} \label{thm:unipotent_Lietype}
Let $p$ be a prime and $G$ a finite simple group of $p$-rank at least two with a strongly $p$-embedded subgroup.
Let $\calC$ be a conjugacy class of non-trivial $p$-elements of $G$.
Then $K_{\calC}$ is irreducible if and only if $\calC$ is not the class of involutions.
\end{theorem}

See Theorem \ref{thm:stronglyPEmbeddedIrreducible}.
Note that, in particular, for conjugacy classes of unipotent elements of order four in the unitary group in characteristic two or in the Suzuki group, the Killing form is irreducible, while for classes of involutions in those same groups it is not.

We address the connectivity of the graph $\killingGraph{G}{\calC}$ by following three different strategies that we outline in Section \ref{sec:overall}.
We then apply these strategies to prove Theorems \ref{thm:involutions_Lietype} and \ref{thm:unipotent_Lietype}.
For the groups of Lie type and Lie rank one, this is done in Sections \ref{sec:linear2} (linear groups), \ref{sec:unitary3} (unitary groups), \ref{sec:suzuki} (Suzuki groups) and \ref{sec:ree} (Ree groups).
In Section \ref{sec:strongly_p_embedded}, we complete the proof of Theorem \ref{thm:unipotent_Lietype} by looking at the remaining groups with a strongly $p$-embedded subgroup from the list in Theorem \ref{thm:stronglypembedded}.
In particular, we need to analyse the case of alternating groups $\Alt_{2p}$ and $\Alt_p$, which is done in Section \ref{sec:symmetric_alternating} (see Theorem \ref{thm:summarySymn_Altn}).
The computations for the remaining groups in the list of Theorem \ref{thm:stronglypembedded} have been done in GAP, in some cases combined with the strategies from Section \ref{sec:overall}.

Finally, in \cite[p.10]{MR3682632}, the authors suggest that every dihedral group of order $2n$, with $n$ odd, is likely to be strongly non-degenerate, and they support this claim with experimental evidence for small values of $n$.
Recall that a group $G$ is strongly non-degenerate if, for every $G$-stable generating subset $\calC\subseteq G\setminus \{1\}$ consisting of real elements, $K_{\calC}$ is non-degenerate.
In Section \ref{sec:dihedral}, we provide a complete proof confirming that dihedral groups of order $2n$, with $n$ odd, are indeed strongly non-degenerate.
See Theorem \ref{thm:dihedralgroup}.

\subsection*{Acknowledgements}
We thank Leandro Vendramin for motivating preliminary discussions that gave rise to this work.

The first author was supported by the FWO grant 12K1223N.
The second author was supported by the FWO Senior Research Project G004124N.
This work was partially supported by the project OZR3762 of Vrije Universiteit Brussel.

\section{Preliminaries}
\label{sec:overall}

In this section, we describe the different strategies that we will follow to show that certain Killing forms on conjugacy classes of finite groups are irreducible.

Let $G$ be a finite group, and let $\calC \subseteq G\setminus \{1\}$ be a $G$-stable subset.
The associated Killing form $K_{\calC}$ is termed \textit{irreducible} if for all $a,b\in \calC$, there is a power $K_{\calC}^k$ of the underlying matrix (written in the basis $\calC$) such that its $(a,b)$-entry is non-zero.
Otherwise, $K_{\calC}$ is termed \textit{reducible}.
In other words, the latter means that there is a suitable ordering of $\calC$ for which $K_{\calC}$ takes the form of a block diagonal matrix.
We will study such a property by looking at the connectivity of a graph.
More precisely, let $\killingGraph{G}{\calC}$ be the graph whose vertex set is $\calC$ and an edge connects two vertices $x,y \in \calC$ if and only if $C_G(xy)\cap \calC\neq \emptyset$, that is, $K_{\calC}(x,y)\neq 0$.
Then, $K_{\calC}$ is irreducible if and only if $\killingGraph{G}{\calC}$ is connected.

To study the connectivity of $\killingGraph{G}{\calC}$, we propose three strategies.

\subsection{Strategy 1}
\label{sub:str1}

The first approach comes from the observation that the commuting graph on $\calC$ is a subgraph of $\killingGraph{G}{\calC}$.

\begin{definition}
Let $G$ be a group, and let $X$ be either a subset of $G$, or a collection of subsets of $G$.
The \textit{commuting graph} on $X$ is the graph $\commGraph{X}$ whose vertex set is $X$, and an edge joins two elements $a,b\in X$ if $[a,b]=1$.
\end{definition}

Note that the commuting graph $\commGraph{\calC}$ is a subgraph of $\killingGraph{G}{\calC}$ that contains every vertex. In particular, the latter is connected if the former is.
However, the converse is not true in general as the commuting graph on $\calC$ might be disconnected while $\killingGraph{G}{\calC}$ is connected.
For example, for $G = \PSL_2(q)$, where $q$ is a power of an odd prime $p$, and $\calC$ a conjugacy class of non-trivial $p$-elements, the commuting graph $\commGraph{\calC}$ is disconnected as a consequence of Theorem \ref{thm:stronglypembedded}, but we will see in Corollary \ref{coro:linear_irred} that $\killingGraph{G}{\calC}$ is connected.

\begin{example} \label{ex:PSL(2,q)_irred}
Let $G$ be a finite simple group of Lie type and Lie rank one in odd characteristic.
Then $G$ has a unique conjugacy class $\calC$ of involutions (see, for instance, \cite{GLS3}), and its commuting graph for this class is connected by Theorem \ref{thm:stronglypembedded}.
Therefore, $K_{\calC}$ is irreducible.
\end{example}

\subsection{Strategy 2}
\label{sub:str2}

We formulate this strategy in an arbitrary $G$-graph, and then we will specialise in the case of graphs $\killingGraph{G}{\calC}$.
Let $\calG$ be a finite simple $G$-graph, where $G$ is a finite group.
Suppose we can cover the vertex set of $\calG$ by subsets $(C_i)_{i\in I}$ (not necessarily disjoint), such that the full subgraph $\calG_i$ spanned by $C_i$ is connected, for $i\in I$ (we shall call these subgraphs \textit{clusters}).
In addition, we will be interested in the case that $G$ transitively permutes the sets $C_i$, hence the subgraphs $\calG_i$.
Next, we form the graph $(\calG, (C_i)_{i\in I})$ with vertex set $I$ and we join $i\neq j\in I$ if either $C_i\cap C_j\neq \emptyset$ or there are vertices $x\in C_i$ and $y\in C_j$ such that $\{x,y\}$ is an edge in $\calG$.
Then, $\calG$ is connected if and only if $(\calG, (C_i)_{i\in I})$ is connected.
If in addition $G$ is double transitive on $\{ C_i : i\in I\}$, then $\calG$ is connected if and only if $(\calG, (C_i)_{i\in I})$ is not a discrete graph with more than one vertex.

Now, suppose we are working with the graph $\killingGraph{G}{\calC}$, for $\calC$ a $G$-conjugacy class.
Let $H$ be a subgroup of $G$ that contains elements of $\calC$, and form the covering $(\calC_{H^g})_{g\in G}$, where $\calC_{H^g} = \calC \cap H^g$, for $g\in G$.
Note that $G$ is transitive on $\{ \calC_{H^g}~:~g\in~G\}$.
Thus, we aim to pick $H$ such that $\calG_H$, the induced graph with vertex set $\calC_H$, and the graph $(\killingGraph{G}{\calC}, (\calC_{H^g})_{g\in G}) $ are connected.
Assume, in addition, that we can take $H$ such that $G$ is double transitive on the $G$-conjugacy class of $H$.
Then, we estimate the number $c_H$ of solutions of $K_{\calC}(x,y)\neq 0$ with $x,y\in \calC_H$.
If, for instance, the sets $\calC_{H^{g}}$ are pairwise disjoint, we know that $c_H \cdot |G:N_G(H)|$ is a lower bound for the number of solutions to the equation $K_{\calC}(x,y)\neq 0$ for $x,y\in \calC$.
Now, one can use counting formulas with characters, such as Theorem \ref{thm:characterFormula}, to show that this bound cannot be attained.
That is, $\killingGraph{G}{\calC}$ cannot be disconnected as otherwise the lower bound $c_H \cdot |G:N_G(H)|$ would be attained.

\begin{theorem}
\label{thm:characterFormula}
Let $G$ be a finite group and let $\calC_1, \calC_2,\calC_3$ be conjugacy classes of $G$.
Then the number of solutions to the equation $xy=z$ with $x\in \calC_1$, $y\in \calC_2$ and $z\in \calC_3$ is
given by 
\begin{equation}
    \Phi_G(\calC_1,\calC_2,\calC_3) := \frac{|\calC_1||\calC_2||\calC_3|}{|G|}
\sum_{\chi\in\Irr(G)}\frac{\chi(\calC_1)\chi(\calC_2)\overline{\chi(\calC_3)}}{\chi(1)}.
\end{equation}
\end{theorem}

Below, we summarise the steps for this strategy.
Let $G$ be a finite group, and $\calC$ a $G$-conjugacy class.
\smallskip

\begin{enumerate}
    \item Pick a suitable subgroup $H\leq G$ such that $\calC\cap H\neq\emptyset$ and $G$ is double transitive on the $G$-conjugacy class of $H$. 
    \item Show that the induced subgraph $\calG_H$ with vertex set $\calC\cap H$ is connected.
    \item Find the $G$-conjugacy classes $\calC_1, \dots, \calC_r$ that intersect $C_G(x)$, for $x \in \calC$.
    \item Prove that there are two distinct conjugates $H^{g_1}$ and $H^{g_2}$ and elements $x\in \calC\cap H^{g_1}$ and $y \in \calC\cap H^{g_2}$ with $K_{\calC}(x,y)\neq 0$ through one of the following methods:
    \begin{enumerate}
        \item Prove that $\calC\cap H^{g_1}\cap H^{g_2}\neq\emptyset$.
        \item By explicit calculations, prove that there exist $x \in \calC \cap H^{g_1}$, $y \in \calC \cap H^{g_2}$ such that $xy \in \calC_i$ for some $i \in \{ 1, \dots, r\}$.
        \item Assume item (a) fails. Then, determine first the number $c_H$ of pairs $(x,y) \in (\calC\cap H)^2$ such that $K_\calC(x,y) \neq 0$. In total, there are thus $c_H \cdot |G : N_G(H)|$ pairs of elements of $\calC$ contained in the same $H^g$ with non-zero Killing form.
        On the other hand, by Theorem \ref{thm:characterFormula}, there are $\sum_{i=1}^r \Phi_G(\calC, \calC, \calC_i)$ pairs $(x,y) \in \calC^2$ such that $K_\calC(x,y) \neq 0$. Compare this number with $c_H  \cdot |G : N_G(H)|$ and conclude that there are pairs of elements contained in distinct $H^g$ with non-zero Killing form.
    \end{enumerate}
    Then, by double transitivity, $\killingGraph{G}{\calC}$ is connected.
\end{enumerate}

\subsection{Strategy 3}
\label{sub:str3}

This strategy relies on computing the stabiliser of a connected component of the graph $\killingGraph{G}{\calC}$.
Let $\calC$ be a non-trivial conjugacy class of $G$.
Since $G$ is transitive on the vertex set of $\killingGraph{G}{\calC}$, $G$ is also transitive on the set of connected components of $\killingGraph{G}{\calC}$.
Let $x\in \calC$, and let $M$ be the stabiliser of the component $C_x$ of $\killingGraph{G}{\calC}$ containing $x$.
Write $\calC_x$ for the vertices contained in $C_x$.
Clearly, $M \geq C_G(y)$ for any vertex $y\in \calC_x$.
More generally, note that if $y\in \calC$ and $g \in N_G(\gen{y})$, then $[y,y^g] = 1$, so $K_{\calC}(y,y^g)\neq 0$ and $y,y^g$ must lie in the same connected component of $\killingGraph{G}{\calC}$.
Hence, we have
\begin{equation}
    \label{eq:componentStab}
    M \geq \langle N_G(\gen{y}) : y \in \calC_x\rangle.
\end{equation}
Moreover, $\killingGraph{G}{\calC}$ is connected if and only if $M = G$.
In particular, if $|\calC_x|> |\calC|/2$, then $M = G$ and $\killingGraph{G}{\calC}$ is connected.

This strategy consists of showing that $M = G$ by using generating properties of $G$, or maximality conditions on the subgroups $N_G(\gen{y})$ for $y\in \calC_x$.

The idea is to apply the strategies described above to $G$, a simple group of Lie type and Lie rank one in characteristic $p$, $H = S$ a Sylow $p$-subgroup of $G$, and $\calC$ a conjugacy class of (non-trivial) unipotent elements, i.e., $p$-elements.
Recall that $G$ is isomorphic to
\[ \PSL_2(p^k), \ \PSU_3(p^k), \ \Sz(2^{2k+1}) \, (p=2), \ \Ree(3^{2k+1}) \, (p=3). \]
In this case, $G$ is double transitive on the set of Sylow $p$-subgroups, which additionally satisfy the trivial intersection property:

\smallskip

\noindent
(TI) \quad Two distinct Sylow $p$-subgroups intersect trivially.

\smallskip

In particular, this implies that $N_G(S)$ is transitive on $\calC\cap S$ for any $S\in \Syl_p(G)$.

The following lemma will allow us to work in the matrix groups $\SL_2(q)$ and $\SU_3(q)$ to check the irreducibility in the projective quotients $\PSL_2(q)$ and $\PSU_3(q)$.

\begin{lemma}
\label{lm:connectedModCenter}
Let $G$ be any group, $Z\leq Z(G)$, and let $\calC$ be a $G/Z$-conjugacy class.
Let $x\in G$ such that the image of $x$ in $G/Z$ lies in $\calC$, and denote by $\calC_x$ its $G$-conjugacy class.
If $\killingGraph{G}{\calC_x}$ is connected, then $\killingGraph{G/Z}{\calC}$ is connected.
Moreover, if the order $o$ of an element of $\calC$ is prime to the order of $Z$, then $x$ can be taken to have order $o$.
\end{lemma}

\begin{proof}
Let $\phi:G\to G/Z$ be the quotient map.
It is clear then that $\phi(\calC_x) = \calC$.
Moreover, if $a,b,c\in \calC_x$ and $c\in C_G(ab)$ then $\phi(c)\in C_{G/Z}(\phi(a)\phi(b))\cap \calC$.
This shows that $\phi$ induces a graph map $\phi_* : \killingGraph{G}{\calC_x} \to \killingGraph{G/Z}{\calC}$ that is surjective on vertices.
The result of the connectivity then follows.

For the Moreover part, suppose $\pi$ is a set of primes such that $o$ is a $\pi$-number and $Z$ is a $\pi'$-group.
Then we have the unique factorisation $x = x_1x_2$ where $x_1$ is the $\pi$-part and $x_2$ is the $\pi'$-part and $x_1,x_2$ commute.
By hypothesis, $x_1^o x_2^o = x^o = z\in Z$, which implies that $x_1^o = 1$ and $x_2^o \in Z$.
Again, since $x_2$ is a $\pi'$-element and $o$ is a $\pi$-number, we can take $o'$ such that $x_2 = (x_2^o)^{o'} = z^{o'} \in Z$.
Thus, $\phi(x_1) = \phi(x)$ and $x_1$ has order exactly $o$.
\end{proof}

We close this section with an example.
As we mentioned in the introduction, Question \ref{ques:reduc} can have a negative answer if the group is not assumed to be simple.
The following example, computed with GAP, shows that we can get reducible Killing forms on conjugacy classes of elements of order three.

\begin{example}
\label{ex:reducible_order3}
Let $G$ be the group with id \texttt{(54,5)} in the library \texttt{SmallGroups} in GAP.
Then, every conjugacy class of elements of order three in $G$ gives rise to a degenerate Killing form, which in some cases is irreducible and in others is not.
On the other hand, $G$ contains a unique conjugacy class of involutions, and the Killing form associated to such a class is reducible and non-degenerate.
\end{example}

\section{Unipotent elements in linear groups of dimension two} \label{sec:linear2}

In this section, we work with the groups $\SL_2(q)$ and $\PSL_2(q)$, where $q$ is a power of a prime $p$.
We take a conjugacy class $\calC$ of non-trivial $p$-elements (i.e., non-trivial unipotent elements), and show that $K_\calC$ is reducible if and only if $p = 2$.

Below we summarise some facts about the linear groups of dimension two.

\begin{lemma} \label{lem:sylowPSL2q}
    Let $q$ be a power of a prime $p$ and $G = \SL_2(q)$ or $\PSL_2(q)$. Then the following hold:
    \begin{enumerate}
        \item $G$ has $q+1$ Sylow $p$-subgroups, each of size $q$.
        \item Each Sylow $p$-subgroup of $G$ is elementary abelian.
        \item Suppose $x$ is a non-trivial unipotent element in $G$ contained in the Sylow $p$-subgroup $S$. Then $C_G(x) = S\times Z(G)$.
        \item If $q$ is odd, then there are exactly two conjugacy classes $\calC_1, \calC_2$ of non-trivial unipotent elements, each of them of size $\frac{q^2 - 1}{2}$ and $|\calC_i\cap S|= \frac{q-1}{2}$ for a Sylow $p$-subgroup $S$.
        Moreover, these classes are real if and only if $q\equiv 1\pmod 4$.
        \item If $q$ is even, there is only one conjugacy class $\calC$ of involutions, which is of size $q^2 - 1$, and thus $|\calC\cap S| = q-1$ for a Sylow $2$-subgroup $S$.
    \end{enumerate}
\end{lemma}

We give some details on these well-known facts, which will serve us to fix the notation.
The set of upper-triangular matrices with ones along the diagonal
\[
S_1 := \biggl\{ \begin{pmatrix}
    1 & a \\
    0 & 1
\end{pmatrix} :  a \in \F_q \biggr\}
\]
is a Sylow $p$-subgroup of $\SL_2(q)$ which projects isomorphically onto a Sylow $p$-subgroup of $\PSL_2(q)$. 

We aim to show that Conjecture \ref{conj:class-non-degen} holds in particular cases, and to answer Question \ref{ques:reduc} in the case of a conjugacy class of unipotent elements in $\PSL_2(q)$.

We start with the case where $q$ is even.
Note that $\SL_2(q) = \PSL_2(q)$ in this case.
We will show that for $\calC$ the conjugacy class of involutions, $K_{\calC}$ is reducible and non-degenerate.
For the latter, we show that certain matrices are invertible, and we will use the following criterion:

\begin{remark}
[{Miller's criterion}]
\label{miller}
Recall from Miller's paper \cite{MR617892} that if $E$ is an invertible matrix and $H$ is a rank-$1$ matrix, then $E+H$ is invertible if and only if $\Tr(HE^{-1})\neq -1$.
In such a case, the inverse of $E+H$ is given by
\[ E^{-1} - \frac{1}{1+g}E^{-1}HE^{-1},\]
\end{remark}

\begin{proposition} \label{prop:psl(2,q)_reduc}
Let $G = \SL_2(q)$, where $q$ is even, and let $\calC$ be the conjugacy class of an involution.
Then $\killingGraph{G}{\calC} = \commGraph{\calC}$, and $K_{\calC}$ is reducible and non-degenerate.
\end{proposition}

\begin{proof}
By Lemma \ref{lem:sylowPSL2q}, there is a single conjugacy class of involutions in $G$, containing therefore all the non-trivial unipotent elements.
In particular, $\calC$ is real and $|\calC| = q^2-1$.
Also, recall that the Sylow $2$-subgroups of $G$ are elementary abelian of order $q$, and they intersect trivially.

Let $x,y,z\in \calC$.
Then $C_G(x)$ equals the unique Sylow $2$-subgroup containing $x$, say $S_x$.
If $yz\in C_G(x) = S_x$ and $y\neq z$, then $yz$ is an involution, that is, $y,z$ commute.
Thus, $y,z$ and $yz$ are contained in the same Sylow $2$-subgroup, namely $S_x$.
This shows that $K_{\calC}(x,y) \neq 0$ if and only if $S_x = S_y$.
In particular, $\killingGraph{G}{\calC} = \commGraph{\calC}$ is disconnected, so $K_\calC$ is irreducible.

Next, we can order $\calC$ according to the Sylow $2$-subgroups.
With such an ordering, the underlying matrix for $K_{\calC}$ is a block diagonal matrix, where the number of blocks is the number of Sylow $2$-subgroups. 
Now, the blocks are all equal to each other and they have dimension $q-1$.
For one such a block, say $B$, we have that the elements on the diagonal are all equal to $q^2-1$ (i.e., the size of $\calC$).
Off the diagonal, we always see the number $| \calC \cap S_x|=q-1$.
Hence, 
\[ B = q(q-1) I_{q-1} + (q-1)\Theta_{q-1,q-1} = (q-1)\big( qI_{q-1} + \Theta_{q-1,q-1}\big),\]
where $\Theta_{m,n}$ denotes the $m\times n$ matrix with all entries equal to one, and $I_m$ is the identity matrix of size $m\times m$.

Let $E = qI_{q-1}$ and $H = \Theta_{q-1,q-1}$.
Note $H$ has rank $1$ and $E$ is invertible.
Also $HE^{-1} = q^{-1}H$ and its trace is $1-q^{-1}\neq -1$.
Thus $B = (q-1)(E+H)$ is invertible by Miller's criterion (see Remark \ref{miller}).
From this, we conclude that $K_{\calC}$ is non-degenerate.
\end{proof}

Now suppose $q$ is odd and $\calC$ is a conjugacy class of non-trivial unipotent elements in $G = \SL_2(q)$ or $\PSL_2(q)$.
We want to show that $K_\calC$ is irreducible by applying Strategy 2 from Subsection \ref{sub:str2}.
Here we take $H$ a Sylow $p$-subgroup of $G$, so the conjugation action on the set of conjugates of $H$ is double transitive.

By Lemma \ref{lem:sylowPSL2q}, the centraliser of an arbitrary element $x \in G$ contains elements from $\calC$ if and only if $x$ is unipotent itself or $x \in Z(G)$.
Hence, for $x,y \in \calC$, the intersection $C_G(xy) \cap \calC$ is non-empty if and only if $xy$ is a product of a unipotent and a central element.
Note that, in the case $G = \SL_2(q)$, $xy\neq -1$.
Next, we analyse the structure of this intersection by using Lemma \ref{lem:sylowPSL2q}.

\smallskip

\noindent
\textbf{Case 1.} $y = x^{-1}$.

Then, we clearly have
\[
K_\calC(x,y) = |C_G(1) \cap \calC| = |\calC| = \frac{q^2 -1}{2}.
\]

\smallskip

\noindent
\textbf{Case 2.} $x$ and $y$ are contained in the same Sylow $p$-subgroup $S$ and $y \neq x^{-1}$.

Then $xy$ is a non-trivial element of $S$. Hence 
\[
K_\calC(x,y) = |C_G(xy) \cap \calC| = |S \cap \calC| = \frac{q-1}{2}.
\]

Cases 1 and 2 imply that any two conjugate elements in the same Sylow $2$-subgroup $S$ have a non-zero Killing form.
In particular, the full subgraph of $\killingGraph{G}{\calC}$ with vertex set $\calC\cap S$ is connected, as we require in the second step of Strategy 2.

\smallskip

\noindent
\textbf{Case 3.} $x,y\in \calC$ lie in distinct Sylow $p$-subgroups.

In this case, a priori, we have no direct conclusion on whether this product lies in the centraliser of an element of $\calC$.
The following proposition determines exactly when such a situation arises.

\begin{proposition} \label{prop:psl(2,q)}
    Let $G = \SL_2(q)$ or $\PSL_2(q)$ with $q$ odd, and let $\calC$ be a conjugacy class of non-trivial unipotent elements in $G$. Let $y_1\in \calC$, and let $S_1$ be the unique Sylow $p$-subgroup containing $y_1$. Then, for each $S_2 \in \Syl_p(G)$ different from $S_1$, there exists a unique $y_2 \in S_2 \cap \calC$ such that $y_1 y_2$ is the product of a non-trivial unipotent element and a central element.
\end{proposition}

\begin{proof}
    Note that $Z(\SL_2(q)) = \{I,-I\}$ has order two, which is prime to the characteristic $p$ of the field $\F_{q}$.
    Suppose that the result holds for $\SL_2(q)$, and denote by $\phi(x)$ the image of an element of $\SL_2(q)$ in $\PSL_2(q)$.
    Take $y_1\in S_1\leq \PSL_2(q)$ as in the statement for $G = \PSL_2(q)$, and fix $S_2$ a Sylow $p$-subgroup of $\PSL_2(q)$, distinct from $S_1$.
    We show that the theorem holds for $\PSL_2(q)$.
    Let $x_1\in \SL_2(q)$ be of the same order as $y_1$ (this is possible by Lemma \ref{lm:connectedModCenter}) and such that $\phi(x_1) = y_1$.
    Note that the quotient map $\SL_2(q) \to \PSL_2(q)$ induces a bijection between the set of Sylow $p$-subgroups, and an isomorphism when restricted to Sylow $p$-subgroups.
    Let $Q_1,Q_2$ be the two Sylow $p$-subgroups of $\SL_2(q)$ that isomorphically map onto $S_1,S_2$ respectively, so $x_1\in Q_1$.
    Then there is $x_2\in Q_2$ such that $x_1x_2$ is the product of a non-trivial unipotent and a central element, so $\phi(x_1x_2)$ is a non-trivial unipotent element, and $y_2 = \phi(x_2) \in S_2$.
    To show the uniqueness, suppose that we have $y_2',y_2\in S_2$ such that $y_1y_2, y_1y_2'$ are non-trivial unipotent.
    Then there exist $x_2,x_2'\in Q_2$ that map onto $y_2,y_2'$ respectively and $(x_1x_2)^p,(x_1x_2')^p$ lie in $Z(\SL_2(q))$.
    That is, both $x_1x_2$ and $x_1x_2'$ can be written as a product of a non-trivial unipotent element and a central element. By uniqueness for $x_1$, we see that $x_2=x_2'$.
    This proves the theorem for $\PSL_2(q)$ if we know it holds for $\SL_2(q)$.
    
    Next, we show that the theorem holds for $\SL_2(q)$.
    By the Sylow Theorem, it is sufficient to prove this claim for 
    \[
    S_1 = \biggl\{ \begin{pmatrix}
        1 & a \\
        0 & 1
    \end{pmatrix} :  a \in \F_q \biggr\}.
    \]

    Recall that $S_1$ acts regularly by conjugation on the set of Sylow $p$-subgroups of $G$ different from $S_1$. By choosing 
    \[
    S_2 := \biggl\{ \begin{pmatrix}
        1 & 0 \\
        x & 1
    \end{pmatrix} :  x \in \F_q \biggr\},
    \]
    as a Sylow $p$-subgroup different from $S_1$, every unipotent element in $G \setminus S_1$ is of the form
    \begin{equation}
    \label{eq:genericUnipotentSL2}
    \begin{pmatrix}
        1 & -a \\
        0 & 1
    \end{pmatrix} \begin{pmatrix}
        1 & 0 \\
        x & 1
    \end{pmatrix} \begin{pmatrix}
        1 & a \\
        0 & 1
    \end{pmatrix} = \begin{pmatrix}
        1 - ax & -a^2x \\
        x & 1 + ax
    \end{pmatrix},
    \end{equation}
    for some $a \in \F_q$ and $x \in \F_q^\times$.

    Let $y_1 = \begin{pmatrix}
        1 & b \\
        0 & 1
    \end{pmatrix} \in \calC$, for $b \in \F_q^\times$, and fix $a \in \F_q$.
    We thus need to prove that there exists a unique $x \in \F_q^\times$ such that
    \[u(x) := \begin{pmatrix}
            1 - ax + bx & -a^2x + b + abx \\
            x & 1 + ax
        \end{pmatrix} =
    \begin{pmatrix}
        1 & b \\
        0 & 1
    \end{pmatrix} \begin{pmatrix}
        1 - ax & -a^2x \\
        x & 1 + ax
    \end{pmatrix}
    \]
    is a product of a non-trivial unipotent and a central element, and such that $y_1$ and
    \[ y_2 = y_2(x) := \begin{pmatrix}
        1 - ax & -a^2x \\
        x & 1 + ax
    \end{pmatrix}\]
    are conjugate. That is, we need to prove that there exist unique $x, z \in \F_q^\times$ and $c \in \F_q$ such that 
    $y_1$ and $y_2$ are conjugate and
    $u(x) =  \pm v(c,z)$, where
    \[ v(c,z) := \begin{pmatrix}
            1 - cz & -c^2z \\
            z & 1 + cz
    \end{pmatrix}. \]
    It is not hard to check that $u(x) = v(c,z)$ has no solutions. On the other hand, $u(x) = -v(c,z)$ has a unique solution given by $x = - \frac{4}{b}$, $z = -x$, and $c = a -\frac{1}{2} b$.
    Furthermore, $y_2(- \frac{4}{b})$ is the conjugate of $y_1$ by
    \[
    \begin{pmatrix}
        0 & \frac{b}{2} \\
        - \frac{2}{b} & 0
    \end{pmatrix}.
    \]
    Thus, the claim holds for $S_1$.
\end{proof}

From this, we can conclude the following result:

\begin{corollary}
\label{coro:linear_irred}
    Let $G = \SL_2(q)$, $\GL_2(q)$, $\PSL_2(q)$ or $\PGL_2(q)$, with $q$ odd, and let $\calC$ be a conjugacy class of non-trivial unipotent elements in $G$. Then $K_\calC$ is irreducible. 
\end{corollary}

\begin{proof}
If $G = \SL_2(q)$ or $\PSL_2(q)$, then the result follows from Proposition \ref{prop:psl(2,q)} and Strategy 2 from Subsection \ref{sub:str2}.
Recall that we take $H =S$ a Sylow $p$-subgroup of $G$, and the subgraph $\calG_S$ is connected since the Sylow $p$-subgroups are abelian (see Strategy 1 in Subsection \ref{sub:str1}).

Suppose now that $G = \GL_2(q)$ or $\PGL_2(q)$.
Here, there is a unique conjugacy class $\calC$ of order-$p$ elements, and the Sylow $p$-subgroups are isomorphic to those of $\SL_2(q)$.
Hence, the same strategy as before shows that they give rise to connected clusters indexed by the Sylow $p$-subgroups.
By Proposition \ref{prop:psl(2,q)}, we can still produce an edge joining two different of these clusters. Thus $\killingGraph{G}{\calC}$ is connected.
\end{proof}

Let $G = \SL_2(q)$ or $\PSL_2(q)$, with $q$ an odd-prime power, and $\calC$ a conjugacy class of order-$p$ elements.
Since the Sylow $p$-subgroups of $G$ intersect trivially, we can again order $\calC$ by the Sylow $p$-subgroups.
Then $K_\calC$ becomes a block matrix $(A_{ij})_{i,j=1,\ldots,q-1}$ such that each $A_{ij}$ is a square matrix of dimension $\frac{q-1}{2}$. By Proposition \ref{prop:psl(2,q)}, for $i \neq j$ we have
\[
A_{ij} = \frac{q-1}{2} P_{ij},
\]
where $P_{ij}$ a permutation matrix. 

The blocks on the diagonal, i.e., when $i = j$, depend on whether or not $\calC$ is a real conjugacy class.
Recall that a $G$-conjugacy class of non-trivial unipotent elements is real if and only if $q \equiv 1 \pmod 4$ (see Lemma \ref{lem:sylowPSL2q}).

Assume first that $q \equiv 1 \pmod 4$. Within each set $\calC\cap S$, where $S$ is a Sylow $p$-subgroup, we can order the elements in a list $x_1,\ldots,x_{2r}$ in such a way that $x_i^{-1}=x_{2r-i+1}$ for all $1\leq i\leq r$. Choosing this ordering on each set $\calC\cap S$, the matrix associated with $K_{\calC}$ has the following blocks on the (main) diagonal:
\[
A_{ii} = \begin{pmatrix}
        \frac{q - 1}{2} & \dots & \frac{q-1}{2} & \frac{q^2-1}{2} \\
        \frac{q-1}{2} & \dots & \frac{q^2-1}{2} & \frac{q-1}{2} \\
        \vdots & \iddots & \vdots & \vdots \\
        \frac{q^2-1}{2} & \dots & \frac{q-1}{2} & \frac{q-1}{2}
    \end{pmatrix}.
\]

Now suppose $q \equiv 3 \pmod 4$. In this case, $\calC$ is not a real conjugacy class, so the product of two conjugate elements in the same Sylow $p$-subgroup can never be the identity. Hence, the blocks on the diagonal are of the form
\[
A_{ii} = \frac{q-1}{2} \Theta_{\frac{q-1}{2}, \frac{q-1}{2}}.
\]

For $q\equiv 1\pmod 4$, one could use the concrete formula we have found above for $K_{\calC}$ and show this is a non-degenerate form, which would confirm Conjecture \ref{conj:class-non-degen} in this particular case.
However, this remains an open problem.

\section{Unipotent elements in special unitary groups of dimension three} \label{sec:unitary3}

In this section, we study Killing forms on conjugacy classes $\calC$ of $p$-elements in $\SU_3(q)$, where $q$ is a power of a prime $p$. Throughout this section, we assume $q \neq 2$.
We show that $K_\calC$ is reducible if and only if $\calC$ is a conjugacy class of involutions.
By Lemma \ref{lm:connectedModCenter}, these results can then be extended to $\PSU_3(q)$.

\subsection{Preliminaries}

The unitary group $\GU_3(q)$ is the subgroup of $\GL_3(q^2)$ consisting of matrices that preserve the form 
\[\Psi((x_1, y_1, z_1), (x_2, y_2, z_2)) = x_1 z_2^q + y_1 y_2^q + z_1 x_2^q.\]
That is, $ \GU_3(q) = \{M \in \GL_3(q^2) : M^* J M = J\},$ where
\[
J = \begin{pmatrix}
    0 & 0 & 1 \\
    0 & 1 & 0 \\
    1 & 0 & 0
\end{pmatrix},
\]
and $M^*$ is the transpose of the matrix obtained by raising each entry of $M$ to the $q$th power.
Then $\SU_3(q)$ is the subgroup of $\GU_3(q)$ of determinant-one matrices.

We recall some general facts about conjugacy classes of $\SU_3(q)$.
We refer to \cite{zbMATH03416061} for more details.
Note that these classes can be computed from the conjugacy classes in $\GU_3(q)$ (see \cite[Section 7]{Ennola}).
However, the canonical forms given in these references are not necessarily unitary matrices.
In Table \ref{tab:PSU(3,q)_conjclasses}, we include suitable unitary representatives of the $\SU_3(q)$-conjugacy classes that we will need.
Remark that the classes $C_6^{(k,l,m)}, C_7^{(k)}$ and $C_8^{(k)}$ are not included in this table, as we do not need further information on them.

Throughout this section, we set $d = \gcd(3,q+1)$, $\xi$ a generator of $\F_{q^2}^\times$, and $\omega = \xi^{\frac{q^2-1}{d}}$ (an element of order $d$).
Recall that $Z(\SU_3(q))$ has order $d$.

\begin{table}[ht]
    \centering
    \begin{tabular}{ | m{1cm} | m{5cm} | m{5.5cm} | }
    \hline
    Class & Representative in $\SU_3(q)$ & Notes \\
    \hline
    $C_1^{(k)}$ & $\begin{pmatrix}
        \omega^k & 0 & 0 \\
        0 & \omega^k & 0 \\
        0 & 0 & \omega^k 
    \end{pmatrix}$ & $0 \leq k \leq d-1$. \\
    $C_2^{(k)}$ & $\begin{pmatrix}
        \omega^k & 0 & \alpha \\
        0 & \omega^k & 0 \\
        0 & 0 & \omega^k 
    \end{pmatrix}$ & $0 \leq k \leq d-1$ and $\alpha \in \F_{q^2}^\times$ such that $\alpha \omega^{kq} + \alpha^q \omega^k = 0$. \\ 
    $C_3^{(k,l)}$ & $\begin{pmatrix}
        \omega^k & \xi^{l(q-1)} & \alpha \\
        0 & \omega^k & - (\xi^l \omega^k)^{-(q-1)}  \\
        0 & 0 & \omega^k 
    \end{pmatrix}$ & $0 \leq k,l \leq d-1$ and $\alpha \in \F_{q^2}^\times$ such that $\omega^{kq}\alpha + \omega^k \alpha^q + 1 = 0$. \\ 
    $C_4^{(k)}$ & $\begin{pmatrix}
        \xi^{k(q-1)} & 0 & 0 \\
        0 & \xi^{-2k(q-1)} & 0 \\
        0 & 0 & \xi^{k(q-1)} 
    \end{pmatrix}$ & $1 \leq k < q+1$, $k \not\equiv 0 \pmod{\frac{q+1}{d}}$. \\
    $C_5^{(k)}$ & $\begin{pmatrix}
        \xi^{k(q-1)} & 0 & \alpha \\
        0 & \xi^{-2k(q-1)} & 0 \\
        0 & 0 & \xi^{k(q-1)} 
    \end{pmatrix}$ & $1 \leq k < q+1$, $k \not\equiv 0 \pmod{\frac{q+1}{d}}$ and $\alpha \in~\F_{q^2}^\times$ satisfies the equation $\alpha \xi^{-k(q-1)} + \alpha^q \xi^{k(q-1)} = 0$. For $q$ even, $\alpha = \xi^{k(q-1)}$ and for $q$ odd, $\alpha = \xi^{k(q-1) - \frac{q+1}{2}}$. \\ 
    \hline
    \end{tabular}
    \caption{Representatives of some of the conjugacy classes of $\SU_3(q)$, $q \neq  2$.}
    \label{tab:PSU(3,q)_conjclasses}
\end{table}

The classes $C_1^{(k)}, 0 \leq k \leq d-1$, constitute the classes of central elements in $\SU_3(q)$.
In particular, if $q \not\equiv -1 \pmod 3$, then $\SU_3(q) \cong \PSU_3(q)$ and the conjugacy classes of these two groups coincide.
If $q \equiv -1 \pmod 3$, the centre of $\SU_3(q)$ is non-trivial, and the $\SU_3(q)$-classes $C_1^{(k)}$, $C_2^{(k)}$ and $C_3^{(k,l)}$ (for a fixed $0\leq l\leq d-1$) collapse into single $\PSU_3(q)$-classes $C_1, C_2$ and $C_3^{(l)}$ respectively.

Table \ref{tab:PSU(3,q)_summary} summarises information on the conjugacy classes of $\SU_3(q)$, including for $C_1^{(0)}, C_2^{(0)}$ and $C_3^{(0,l)}$ the number of elements of the class contained in a given Sylow $p$-subgroup. In the last column, we include the conjugacy classes that intersect the centraliser of an element in the class of the row. We extract this information from Lemmas \ref{lem:PSU_3q_cent} and \ref{lem:PSU_3q_centC3} below.
An empty cell in the table means that the information is not needed and thus not included.

\begin{table}[ht]
    \centering
    \begin{tabular}{| m{1.15cm} | m{2.5cm} | m{2.6cm} | m{1.25cm} | m{3cm} |}
        \hline
       Class & $|x|$ for $q$ even/odd & $|C_G(x)|$ & $|\calC\cap S|$ & Classes in $C_G(x)$ \\
       \hline
       $C_1^{(k)}$ & $1$ if $k=0$, and $3$ otherwise & $q^3 (q^3 +1) (q^2 -1)$ & $1$ & all \\
       $C_2^{(k)}$  &  $2/p$ if $k=0$, $6/3p$ otherwise & $q^3 (q+1)$ & $q-1$ & $C_1^{(k')},C_2^{(k')}, C_3^{(k',l)}$, $C_4^{(k')}, C_5^{(k')}$ \\
       $C_3^{(k,l)}$ & $4/p$ if $k=0$, $12/3p$ otherwise & $dq^2$ & $q(q^2-1)$ &  $C_1^{(k')},C_2^{(k')},C_3^{(k',l)}$ \\
       $C_4^{(k)}$ & divides $q+1$ & $q(q+1)^2(q-1)$ &  &  \\
       $C_5^{(k)}$ & divides $p(q+1)$ & $q(q+1)$ &  &  \\
       $C_6^{(k,l,m)}$ & divides $q+1$ & $(q+1)^2$ &  &  \\
       $C_7^{(k)}$ & divides $q^2-1$ & $q^2-1$ &  &  \\
       $C_8^{(k)}$ & divides $q^2-q+1$ & $q^2-q+1$ &  & \\
       \hline
    \end{tabular}
    \caption{List of conjugacy classes $\calC$ in $\SU_3(q)$, $q \neq 2$. Here, $S$ is a Sylow $p$-subgroup of $\SU_3(q)$, and $x$ is a representative of the corresponding class.}
    \label{tab:PSU(3,q)_summary}
\end{table}

Note that $C_2^{(0)}$ and $C_3^{(0,l)}$, $0 \leq l \leq d-1$, are the conjugacy classes of unipotent elements in $\SU_3(q)$. Thus, they make up the non-trivial elements of the Sylow $p$-subgroups, where the elements in the centres are given by $C_2^{(0)}$. 
For $d = 3$, the conjugacy classes of unipotent elements $C_3^{(0,l)}, 0 \leq l \leq 2$, stem from a single conjugacy class in $\GU_3(q)$, as the following lemma explains.

\begin{lemma} \label{lem:SU3(q)_C3_conjugates}
    Let $G = \GU_3(q)$ with $q \equiv -1 \pmod 3$. Then there exist $g_l \in G$ such that $C_3^{(0,l)} = g_l  C_3^{(0,0)} g_l^{-1}$ for $l=1,2$.
\end{lemma}

\begin{proof}
    If $\xi$ is a generator of $\F_{q^2}^\times$, then take
    \[
    g_1 = \begin{pmatrix}
        1 & 0 & 0 \\
        0 & \xi^{1-q} & 0 \\
        0 & 0 & 1
    \end{pmatrix} \qquad \text{and} \qquad g_2 = \begin{pmatrix}
        1 & 0 & 0 \\
        0 & \xi^{2 - 2q} & 0 \\
        0 & 0 & 1
    \end{pmatrix}.
    \]
\end{proof}

Both groups $\SU_3(q)$ and $\PSU_3(q)$ have $q^3+1$ Sylow $p$-subgroups, each of size $q^3$. 
As in the linear case, with our choice of the Hermitian form, the set of
upper-triangular unitary matrices with ones on the diagonal constitute a Sylow $p$-subgroup of $\SU_3(q)$, which we denote by $S_1$:
\[
S_1 := \Biggl\{ \begin{pmatrix}
    1 & a_1 & a_2 \\
    0 & 1 & a_3 \\
    0 & 0 & 1
\end{pmatrix} : \; a_1, a_2, a_3 \in \F_{q^2} \text{ such that } a_3 = -a_1^q, \; a_2 + a_2^q  + a_1^{q+1} = 0 \Biggr\}.
\]
Note that $S_1$ is also a Sylow $p$-subgroup of $\GU_3(q)$ which projects isomorphically onto a Sylow $p$-subgroup of $\PSU_3(q)$.
The centre of $S_1$ consists of the matrices for which $a_1 = a_3 = 0$ and $a_2+a_2^q = 0$.
Moreover, $S_1$ acts regularly by conjugation on the set of Sylow $p$-subgroups different from $S_1$ (recall that $S_1$ is the unipotent radical of its normaliser, which is a Borel subgroup), and any of the groups $\SU_3(q)$, $\PSU_3(q)$, $\GU_3(q)$ and $\PGU_3(q)$ satisfies (TI).

Similarly, we will work with $S_2$, the Sylow $p$-subgroup spanned by the lower-triangular unitary matrices with ones on the diagonal:
\[
S_2 := \Biggl\{ \begin{pmatrix}
    1 & 0 & 0 \\
    b_1 & 1 & 0 \\
    b_2 & b_3 & 1
\end{pmatrix} : \; b_1, b_2, b_3 \in \F_{q^2} \text{ such that } b_3 = -b_1^q, \; b_2 + b_2^q + b_1^{q+1} = 0 \Biggr\}.
\]

We aim to apply Strategies 2 and 3 from Subsections \ref{sub:str2} and \ref{sub:str3} respectively, to show that the Killing form on a conjugacy class of non-involutive unipotent elements is irreducible.
Considering Strategy 3, we recall the following result from \cite[Theorem 3]{Nuzhin} for $\PSU_3(q) = \SU_3(q)/Z(\SU_3(q))$.
Let $\phi:\SU_3(q) \rightarrow~\PSU_3(q)$ denote the natural quotient map.

\begin{theorem} \label{thm:Nuzhin_gen_PSU3(q)}
    For $q > 4$, $\PSU_3(q)$ is generated by $\phi(u_0)$ and $\phi(Z(S_2))$, where
    \[
    u_0 = \begin{pmatrix}
        1 & 1 & \alpha \\
        0 & 1 & -1 \\
        0 & 0 & 1
    \end{pmatrix},
    \]
    with $\alpha\in \F_{q^2}^\times$ satisfying $\alpha + \alpha^q + 1 = 0$.
    
    In particular, $\langle u_0, Z(S_2), Z(\SU_3(q)) \rangle = \SU_3(q)$.
\end{theorem}

Note that $u_0$ is a representative of $C_3^{(0,0)}$ as given in Table \ref{tab:PSU(3,q)_conjclasses}.

\subsection{The class \texorpdfstring{$C_2^{(0)}$}{C\_2\^{(0)}}}

We first consider Killing forms on the conjugacy class of order-$p$ elements, namely, $\calC = C_2^{(0)}$. 
Following Strategy 2, we study the centralisers of elements in $\calC$.

\begin{lemma} \label{lem:PSU_3q_cent}
    Let $G = \SU_3(q)$ and $x \in C_2^{(0)}$. Then $C_G(x)$ only contains elements from the conjugacy classes $C_1^{(k)}, C_2^{(k)}, C_3^{(k,l)}, C_4^{(k)}$ and $C_5^{(k)}$.
\end{lemma}

\begin{proof}
    It is enough to prove the claim for $x$ the representative of $C_2^{(0)}$ given in Table \ref{tab:PSU(3,q)_conjclasses}, with $\alpha+\alpha^q=0$, $\alpha\in \F_{q^2}^{\times}$. Then, an element 
    \[ y = \begin{pmatrix}
    y_{11} & y_{12} & y_{13} \\
    y_{21} & y_{22} & y_{23} \\
    y_{31} & y_{32} & y_{33}
    \end{pmatrix} \in \SU_3(q)\]
    commutes with $x$ if and only if
    \[
    \begin{pmatrix}
        y_{11} + \alpha y_{31} & y_{12} + \alpha y_{32} & y_{13} + \alpha y_{33} \\
        y_{21} & y_{22} & y_{23} \\
        y_{31} & y_{32} & y_{33}
    \end{pmatrix} = \begin{pmatrix}
        y_{11} & y_{12} & \alpha y_{11} + y_{13} \\
        y_{21} & y_{22} & \alpha y_{21} + y_{23} \\
        y_{31} & y_{32} & \alpha y_{31} + y_{33}
    \end{pmatrix}.
    \]
    This forces $y_{21}=0=y_{31}=y_{32}$ and $y_{11} = y_{33}$.
    Furthermore, as $1 = \det(y) = y_{11}^2 y_{22}$, we get
    \[
    y = \begin{pmatrix}
        y_{11} & y_{12} & y_{13} \\
        0 & y_{11}^{-2} & y_{23} \\
        0 & 0 & y_{11}
    \end{pmatrix}.
    \]
    The eigenvalues of $y$ are thus $y_{11}, y_{11}, y_{11}^{-2}$, so every element in the centraliser of an involution must have a spectrum of the form $\{a, a, a^{-2}\}$, with $a\in \F_{q^2}^\times$.
    Since the spectrum is fixed under conjugation and the conjugacy classes of $\SU_3(q)$ arise from those of $\GU_3(q)$, we can compare this to the spectra of the canonical representations of the classes given in \cite[Section 7]{Ennola}. 
    From this, we conclude that the centraliser of an element in $C_2^{(0)}$ only contains elements from the conjugacy classes $C_1^{(k)}, C_2^{(k)}, C_3^{(k,l)}, C_4^{(k)}$ and $C_5^{(k)}$.
\end{proof}

Now we focus on the class $C_2^{(0)}$ when $q$ is even. 
Recall from Table \ref{tab:PSU(3,q)_summary} that $\calC =  C_2^{(0)}$ is then the unique class of involutions in $\SU_3(q)$. We show that $K_{\calC}$ is reducible and non-degenerate. 

\begin{proposition} \label{prop:PSU(3,q)_reduc}
Let $G = \SU_3(q)$ or $\PSU_3(q)$ with $q$ even, and let $\calC$ be the conjugacy class of involutions in $G$. Then $\killingGraph{G}{\calC} = \commGraph{\calC}$, and $K_\calC$ is reducible and non-degenerate.
\end{proposition}

\begin{proof}
We assume that $G = \SU_3(q)$, so $\calC = C_2^{(0)}$.
By Lemma \ref{lem:PSU_3q_cent}, $C_G(x) \cap \calC$ is non-empty only if $x$ lies in one of the conjugacy classes $C_1^{(k)}, C_2^{(k)}, C_3^{(k,l)}, C_4^{(k)}$ and $C_5^{(k)}$.

We look at the conjugacy classes which can be represented by elements that are products of two involutions.
Let $x,y \in \calC$. Then $xy = y(yx)y$ is real, so by \cite[Lemma 2]{zbMATH03197179},
$xy$ is either a $2$-element, thus contained in classes $C_1^{(0)}, C_2^{(0)}$ or $C_3^{(0,k)}$, or it has a centralizer of odd order, thus contained in $C_6^{(k,l,m)}, C_7^{(k)}$ or $C_8^{(k)}$ (see Table \ref{tab:PSU(3,q)_summary}).
Therefore, for $x,y \in \calC$, $K_\calC(x,y) \neq 0$ if and only if $xy$ is a $2$-element, and so $xy \in C_1^{(0)}, C_2^{(0)}$ or $C_3^{(0,k)}$.

Suppose $xy$ is a $2$-element. Since $\gen{x,y}$ is a $2$-group and contains both $x,y$, by (TI) it must lie in a unique Sylow $2$-subgroup, say $S$. Hence, $x,y$ are involutions in the centre of $S$, implying that $xy$ has order at most two.
This proves that $xy$ is a $2$-element if and only if it belongs to $C_1^{(0)}$ or $\calC = C_2^{(0)}$, if and only if $x,y$ belong to the same Sylow $2$-subgroup.

We conclude that for $x, y \in \calC$, $K_\calC(x,y) = |C_G(xy) \cap \calC| \neq 0$ if and only if $x,y$ lie in the same Sylow $2$-subgroup. 
We then have $\killingGraph{G}{\calC} = \commGraph{\calC}$ and
\[ K_\calC(x,y) = \begin{cases}
   |\calC| & x = y,\\ 
   |\calC\cap S| & x\neq y, x,y\in S\in \Syl_2(G),\\
   0 & \text{otherwise.}
\end{cases}
\]
By Table \ref{tab:PSU(3,q)_summary}, $|\calC| = q^3(q-1)$ and $|S \cap \calC| = q-1$.

We can order the elements in $\calC$ by the Sylow $2$-subgroups they are contained in. 
The matrix of the Killing form on $\calC$ then becomes a block diagonal matrix.
The blocks on the diagonal are all equal, of dimension $q-1$, and with entries
\[
\begin{pmatrix}
    (q^3+1)(q-1) & q-1 & \dots & q-1 \\
    q-1 & (q^3+1)(q-1)  & \dots & q-1 \\
    \vdots & \vdots & \ddots & \vdots \\
    q-1 & q-1 & \dots & (q^3+1)(q-1) 
\end{pmatrix}.
\]

By Miller's criterion (see Remark \ref{miller}), with $E = (q^3+1)(q-1) I_{q-1}$ and $H = (q-1)\Theta_{q-1,q-1}$, these matrices on the diagonal are invertible. 
Therefore, $K_\calC$ is reducible and non-degenerate.

The result for $G = \PSU_3(q)$ follows similarly by noting that $C_2^{(0)}$ projects isomorphically to the $\PSU_3(q)$-conjugacy class of involutions, and for an involution $x\in \SU_3(q)$, $C_{\PSU_3(q)}(\phi(x)) = \phi( C_{\SU_3(q)}(x))$, where $\phi:\SU_3(q)\to \PSU_3(q)$ is the quotient map.
Moreover, $\phi$ induces a graph isomorphism $\killingGraph{\SU_3(q)}{C_2^{(0)}} \cong \killingGraph{\PSU_3(q)}{\phi(C_2^{(0)})}$.
\end{proof}

Next, assume $q$ is an odd-prime power and let again $\calC = C_2^{(0)}$.
Suppose $x,y \in \calC$. By Lemma \ref{lem:PSU_3q_cent}, $K_\calC(x,y) = |C_G(xy) \cap \calC| \neq 0$ if and only if $xy$ lies in $C_1^{(k)}, C_2^{(k)}, C_3^{(k,l)}$, $C_4^{(k)}$ or $C_5^{(k)}$. 
If $x$ and $y$ are contained in the same Sylow $p$-subgroup $S$, $x,y \in Z(S)$ and thus $xy \in Z(S)$. Hence, in this case, $xy \in C_2^{(0)}$ or $xy=1$. Therefore, for all $x,y \in \calC$ that lie in the same Sylow $p$-subgroup, we have $K_\calC(x,y) \neq 0$. 

Following Strategy 2, we now check when two elements from distinct Sylow $p$-subgroups have non-zero Killing form.
We have the following proposition.

\begin{proposition} \label{prop:PSU3_C2odd}
    Let $G = \SU_3(q)$ for $q$ odd, and let $x \in C_2^{(0)}$. Suppose $x$ lies in the Sylow $p$-subgroup $S$ and let $S' \in \Syl_p(G) \setminus \{S\}$.
    Take $y \in S' \cap C_2^{(0)}$.
    Then $K_{\calC}(x,y) \neq 0$ if and only if $xy \in C_5^{(k)}$ with $k = \frac{q+1}{2}$, and for any $x$, such an element $y$ always exists.
\end{proposition}
 
\begin{proof}
By the Sylow Theorem, we may assume that $S = S_1$ consists of upper-triangular unitary matrices. Recall that $S_2$ denotes the Sylow $p$-subgroups consisting of lower-triangular unitary matrices.
We will prove the claim for $x \in S_1$ and $S' = S_2$. Because of the regularity of the action of $S_1$ on $\Syl_p(G) \setminus \{S_1\}$, the conclusions of the proposition will then hold true for every Sylow $p$-subgroup different from $S_1$. 

Let $y \in S_2 \cap C_2^{(0)}$. As $x$ and $y$ are central in $S_1$ and $S_2$ respectively, they are of the form
\[
x = \begin{pmatrix}
    1 & 0 & x_2 \\
    0 & 1 & 0 \\
    0 & 0 & 1
\end{pmatrix}, \quad y = \begin{pmatrix}
    1 & 0 & 0 \\
    0 & 1 & 0 \\
    y_2 & 0 & 1
\end{pmatrix},
\]
for some $x_2, y_2 \in \F_{q^2}^\times$ with $x_2+x_2^q=0$ and $y_2+y_2^q =0$. Hence 
\[
xy = \begin{pmatrix}
    1 + x_2 y_2 & 0 & x_2 \\
    0 & 1 & 0 \\
    y_2 & 0 & 1
\end{pmatrix}.
\]
Remark that as $x_2, y_2 \neq 0$, one can verify that the minimal polynomial of $xy$ must have degree three (so equal to the characteristic polynomial),
with $1$ as a root of multiplicity one and two other roots different from $1$ (but possibly equal to each other).
By Lemma \ref{lem:PSU_3q_cent}, if $K_{\calC}(x,y)\neq 0$, then $xy$
lies in $C_1^{(k)}, C_2^{(k)}$, $C_3^{(k,l)}$, $C_4^{(k)}$ or $C_5^{(k)}$.
We show that the only possible case is $xy \in C_5^{(k)}$ with $k = \frac{q+1}{2}$, and that such an element $y$ with $K_{\calC}(x,y)\neq 0$ exists.

To see this, note that representatives of the classes $C_1^{(k)}, C_2^{(k)}$ and $C_3^{(k,l)}$ have a unique eigenvalue of multiplicity three, which is not the case for $xy$.
Also, representatives in $C_4^{(k)}$ for $0 \leq k < q+1$, have a minimal polynomial of degree two, in contrast to the degree-three minimal polynomial of $xy$.
Therefore, if $K_{\calC}(x,y)\neq 0$ then $xy$ must lie in $C_5^{(k)}$.

Suppose then that $xy \in C_5^{(k)}$ for some $0 \leq k < q+1$.
As $1$ is an eigenvalue of $xy$, we must have $k = \frac{q+1}{2}$.
Hence, $xy = a^{-1} b a$, where
\[
a = \begin{pmatrix} 
    a_{11} & a_{12} & a_{13} \\
    a_{21} & a_{22} & a_{23} \\
    a_{31} & a_{32} & a_{33}
\end{pmatrix} \in \SU_3(q), \quad \text{ and } \quad b = \begin{pmatrix}
    -1 & 0 & \alpha \\
    0 & 1 & 0 \\
    0 & 0 & -1
\end{pmatrix},\]
with $\alpha = -\xi^{- \frac{q+1}{2}}$ and $a_{ij} \in \F_{q^2}, \; 1 \leq i,j \leq 3$. Thus, from $axy=ba$ we get to

\begin{multline*}
    \begin{pmatrix}
    a_{11} (1+x_2 y_2) + a_{13} y_2 & a_{12} & a_{11} x_2 + a_{13} \\
    a_{21} (1+x_2 y_2) + a_{23} y_2 & a_{22} & a_{21} x_2 + a_{23} \\
    a_{31} (1+x_2 y_2) + a_{33} y_2 & a_{32} & a_{31} x_2 + a_{33}
\end{pmatrix} \\
= \begin{pmatrix}
    -a_{11} + \alpha a_{31}  & - a_{12} + \alpha a_{32} & -a_{13} + \alpha a_{33} \\
    a_{21} & a_{22} & a_{23} \\
    -a_{31} & -a_{32} & -a_{33}
\end{pmatrix}.
\end{multline*}

The only solution to this equation is $y_2 = - \frac{4}{x_2}$, $a_{12} = a_{21} = a_{23} = a_{32} = 0$, $a_{13} = \frac{a_{33} \alpha - a_{11} x_2}{2}$, and $a_{31} = - \frac{2 a_{33}}{x_2}$. 
As $1 = \det(a) = a_{22}(a_{11}a_{33} - a_{13}a_{31})$, we get $a_{22} = \frac{x_2}{\alpha a_{33}^2}$.
Furthermore, for $a$ to be unitary, the following should hold:
\[
\begin{pmatrix}
    a_{11}^q & 0 & a_{31}^q \\
    0 & a_{22}^q & 0 \\
    a_{13}^q & 0 & a_{33}^q
\end{pmatrix} \begin{pmatrix}
    0 & 0 & 1 \\
    0 & 1 & 0 \\
    1 & 0 & 0
\end{pmatrix} \begin{pmatrix}
    a_{11} & 0 & a_{13} \\
    0 & a_{22} & 0 \\
    a_{31} & 0 & a_{33}
\end{pmatrix} = \begin{pmatrix}
    0 & 0 & 1 \\
    0 & 1 & 0 \\
    1 & 0 & 0
\end{pmatrix}.
\]
This further implies that $a_{33}^{q+1} = \frac{x_2}{\alpha}$ and $a_{11}^{q-1} = a_{33}^{q-1}$.
Notice that $(\frac{x_2}{\alpha})^{q-1} = 1$, so we can find $a_{33}$ satisfying $a_{33}^{q+1} = \frac{x_2}{\alpha}$, and take $a_{11} = a_{33}$.
Since the other coefficients of $a$ are determined from $a_{11},a_{33}$, such a matrix $a$ exists.

This shows that there exists $y \in S_2 \cap C_2^{(0)}$ such that $xy \in C_5^{(k)}$ and $K_{\calC}(x,y)\neq 0$, if and only if $k = \frac{q+1}{2}$.
\end{proof}

From Proposition \ref{prop:PSU3_C2odd} and Strategy 2 from Subsection \ref{sub:str2}, we conclude:

\begin{corollary}
\label{coro:psu_C2_qodd_irred}
Let $G = \SU_3(q)$ with $q$ odd, and let $\calC = C_2^{(0)}$.
Then $K_{\calC}$ is irreducible.
\end{corollary}

\begin{proof}
We follow Strategy 2 using a Sylow $p$-subgroup $S$ as our group $H$ there.

Since $C_2^{(0)} \cap S$ consists of central elements in $S$, the clusters generated by these sets are connected subgraphs.

On the other hand, from Proposition \ref{prop:PSU3_C2odd}, for all $x \in C_2^{(0)} \cap S$, and for each $S' \in \Syl_p(G) \setminus \{S\}$, there exists $y \in S'$ such that $K_{C_2}(x,y) \neq 0$.
Therefore, the graph $\killingGraph{G}{\calC}$ is connected.
\end{proof}

\subsection{The classes \texorpdfstring{$C_3^{(0,l)}$ }{C\_3\^{(0,l)}}}

In this subsection, we prove that the Killing form on the conjugacy classes $C_3^{(0,l)} \; (0 \leq l \leq d-1)$ is irreducible for any $q$. 
In the case $q$ is even, $C_3^{(0,l)}$ are conjugacy classes of elements of order four.

We first require the following lemma.

\begin{lemma} \label{lem:PSU(3,q)_C3conjugates}
    Suppose $q \equiv -1 \pmod 3$. Let
    \[
    x = \begin{pmatrix}
        1 & x_1 & x_2 \\
        0 & 1 & x_3 \\
        0 & 0 & 1
    \end{pmatrix}, \quad y = \begin{pmatrix}
        1 & y_1 & y_2 \\
        0 & 1 & y_3 \\
        0 & 0 & 1
    \end{pmatrix}
    \]
    be elements in $S_1 \leq \SU_3(q)$ and $\xi$ a generator of $\F_{q^2}^\times$. Then $x$ and $y$ are conjugate if and only if $x_1 = \xi^{3m} y_1$ for some $0 \leq m < \frac{q^2-1}{3}$.
\end{lemma}

\begin{proof}
    The normaliser $N_G(S_1)$ of $S_1$ in $G = \SU_3(q)$ consists of elements of the form
    \[
    \begin{pmatrix}
        a_{11} & a_{12} & a_{13} \\
        0 & a_{22} & a_{23} \\
        0 & 0 & a_{33}
    \end{pmatrix}.
    \]
    In order for them to lie in $\SU_3(q)$, the coefficients of these matrices need to satisfy, among others, the additional conditions
    \[
    a_{11} a_{22} a_{33} = 1, \quad \text{ and } \quad a_{11}^q a_{33} = 1.
    \]
    In particular, this implies that $a_{22} = (a_{11} a_{33})^{-1} = a_{11}^{q-1}$.

    Therefore, conjugates of $x$ within $S_1$ are of the form
    \begin{multline*}
    \begin{pmatrix}
        a_{11} & a_{12} & a_{13} \\
        0 & a_{22} & a_{23} \\
        0 & 0 & a_{33}
    \end{pmatrix} ^{-1} \begin{pmatrix}
        1 & x_1 & x_2 \\
        0 & 1 & x_3 \\
        0 & 0 & 1
    \end{pmatrix}
    \begin{pmatrix}
        a_{11} & a_{12} & a_{13} \\
        0 & a_{22} & a_{23} \\
        0 & 0 & a_{33}
    \end{pmatrix} \\
    = \begin{pmatrix}
        1 & x_1 a_{22}^2 a_{33} & x_1 a_{22} a_{23} a_{33} + x_2 a_{22}^2 a_{33} - x_3 a_{12}a_{33}^2 \\
        0 & 1 & x_3 a_{11} a_{33}^2 \\
        0 & 0 & 1
    \end{pmatrix}.
    \end{multline*}

    Thus 
    \[
    \frac{x_1 a_{22}^2 a_{33}}{x_1} = a_{22}^2 a_{33} = (a_{11} a_{22} a_{33}) a_{11}^{-1} a_{22} = a_{11}^{q-2}.
    \]
    As $q \equiv -1 \pmod 3$, $q-2$ is divisible by three, so $a_{11}^{q-2}$ is a cube. This proves one implication of the claim. 

    On the other hand, one third of the elements in $S_1 \setminus Z(S_1)$ are conjugate to $x$, and one third of the elements in $S_1 \setminus Z(S_1)$ satisfy the condition from the claim. This proves that every element of $S_1$ conjugate to $x$ must be of this form.
\end{proof}

Recall that the conjugacy classes $C_3^{(k,l)}$, for a fixed $l$, collapse to a single conjugacy class $C_3^{(l)}$ in $\PSU_3(q)$. 
Therefore, we can use Lemma \ref{lem:PSU(3,q)_C3conjugates} to verify whether, for $x \in C_3^{(0,l)}$ and $y \in C_3^{(k,l')}$, their images in $\PSU_3(q)$ are conjugate.
Recall that $\phi: \SU_3(q) \rightarrow \PSU_3(q)$ is the canonical quotient map.

\begin{corollary} \label{cor:SU3_C3kl_samel}
    Suppose $q \equiv -1 \pmod 3$. Let
    \[
    x = \begin{pmatrix}
        1 & x_1 & x_2 \\
        0 & 1 & x_3 \\
        0 & 0 & 1
    \end{pmatrix} \in C_3^{(0,l)}, \quad \text{and} \quad y = \begin{pmatrix}
        \omega^k & y_1 & y_2 \\
        0 & \omega^k & y_3 \\
        0 & 0 & \omega^k
    \end{pmatrix} \in C_3^{(k,l')}.
    \]
    Then the following are equivalent:
    \begin{enumerate}
        \item $y_1 = \xi^{3m} \omega^k x_1$ for some $0 \leq m < \frac{q^2-1}{3}$;
        \item $\phi(x)$ and $\phi(y)$ are conjugate in $\PSU_3(q)$;
        \item $l = l'$.
    \end{enumerate}
\end{corollary}

Following Strategy 2, in order to show that the Killing form on $C_3^{(0,l)}$ is irreducible, we use Lemma \ref{lem:PSU(3,q)_C3conjugates} and Corollary \ref{cor:SU3_C3kl_samel} to study the centralisers of elements in $C_3^{(0,l)}$.

\begin{lemma} \label{lem:PSU_3q_centC3}
    Let $G=\SU_3(q)$, $0 \leq l \leq d-1$, and $x \in C_3^{(0,l)}$. Then $C_G(x)$ only contains elements from $C_1^{(k)}, C_2^{(k)}$ and $C_3^{(k,l)}$, for $0\leq k\leq d-1$.
\end{lemma}

\begin{proof}
We prove the claim for 
\[
x = \begin{pmatrix}
    1 & x_1 & x_2 \\
    0 & 1 & x_3 \\
    0 & 0 & 1
\end{pmatrix} \in S_1 \cap C_3^{(0,l)}.
\]
Then $y = (y_{ij}) \in \SU_3(q)$ commutes with $x$ if and only if
\begin{multline*}
    \begin{pmatrix}
        y_{11} + x_1 y_{21} + x_2 y_{31} & y_{12} + x_1 y_{22} + x_2 y_{32} & y_{13} + x_1 y_{23} + x_2 y_{33} \\
        y_{21} + x_3 y_{31} & y_{22} + x_3 y_{32} & y_{23} + x_3 y_{33} \\
        y_{31} & y_{32} & y_{33}
    \end{pmatrix} \\ = \begin{pmatrix}
        y_{11} & x_1 y_{11} + y_{12} & x_2 y_{11} + x_3 y_{12} + y_{13} \\
        y_{21} & x_1 y_{21} + y_{22} & x_2 y_{21} + x_3 y_{22} + y_{23} \\
        y_{31} & x_1 y_{31} + y_{32} & x_2 y_{31} + x_3 y_{32} + y_{33}
    \end{pmatrix}.   
\end{multline*}
Solving the corresponding system of equations yields that $y$ must be of the form
\[
y = \begin{pmatrix}
    y_{11} & y_{12} & y_{13} \\
    0 & y_{11} & \frac{x_3}{x_1} y_{12} \\
    0 & 0 & y_{11}
\end{pmatrix} = \begin{pmatrix}
    \omega^k & y_1 & y_2 \\
    0 & \omega^k & y_3 \\
    0 & 0 & \omega^k
\end{pmatrix},
\]
with $\omega^k = y_{11}$, $y_1 = y_{12}$, $y_2 = y_{13}$ and $y_3 = \frac{x_3}{x_1} y_1$.
If $y_1=y_2=y_3=0$, $y \in C_1^{(k)}$. If $y_1=y_3 = 0$ and $y_2 \neq 0$, $y \in C_2^{(k)}$. 
If the $y_i$ are all non-zero, $y \in C_3^{(k,l')}$ for some $0 \leq l'\leq d-1$.

It remains to prove that if $y\in C_3^{(k,l')}$, then $l' = l$.
If $d=1$, this is clear. 
Suppose thus that $d=3$, that is, $q\equiv -1\pmod{3}$. 
As $x$ and $y$ are both unitary, $x_3 = -x_1^q$ and $y_3 = - \omega^{k(1-q)} y_1^q$.
Therefore, 
\[\frac{y_1}{x_1} = \frac{y_3}{x_3} = \frac{- \omega^{k(1-q)} y_1^q}{-x_1^q} = \Big( \frac{y_1}{x_1} \Big)^q \omega^{k(1-q)}.
\]
Dividing both sides by $\frac{y_1}{x_1}$, we see 
\[
1 = \Big( \frac{y_1}{x_1} \Big)^{q-1} \omega^{k(1-q)} \Longleftrightarrow \Big( \frac{y_1}{x_1} \Big)^{q-1} = \omega^{k(q-1)},
\]
or thus $\frac{y_1}{x_1} = \omega^k \xi^{m(q+1)}$ for some $0 \leq m < q-1$.
As $3$ divides $q+1$, by Corollary \ref{cor:SU3_C3kl_samel}, $l = l'$.
\end{proof}

Note that Lemma \ref{lem:PSU(3,q)_C3conjugates}, Corollary \ref{cor:SU3_C3kl_samel} and Lemma \ref{lem:PSU_3q_centC3} hold for any $q$ regardless of its parity.

By Lemma \ref{lem:PSU_3q_centC3}, for $x,y \in \calC = C_3^{(0,l)}$ with $l$ fixed, $K_\calC(x,y) \neq 0$ if and only if $xy \in C_1^{(k)}, C_2^{(k)}$ or $C_3^{(k,l)}$ for some $k$.
We will see that, when $q\not\equiv -1\pmod 3$, we can follow Strategy 2 from Subsection \ref{sub:str2} to show that $K_{\calC}$ is irreducible for $\calC = C_3^{(0,0)}$.
For $q \equiv -1 \pmod 3$, we will have to apply Strategy 3 from Subsection \ref{sub:str3}.
Either way, we need to show that there are connections between vertices that lie in different Sylow $p$-subgroups in the graph $\killingGraph{\SU_3(q)}{\calC}$.
To this end, we have the following proposition.

\begin{proposition} \label{prop:SU3_C3}
    Let $G = \SU_3(q)$ with $q\neq 2$, $\calC = C_3^{(0,0)}$ and $S, S'$ two distinct Sylow $p$-subgroups of $G$. Then there are $x \in S \cap \calC$ and $y \in S' \cap \calC$ such that $xy \in \calC$.

    Moreover, when $S = S_1$ and $S' = S_2$, we can take $x\in S_1$ to be as follows:
    \[ x = \begin{cases}
    \begin{pmatrix}
    1 & 1 & x_2 \\
    0 & 1 & -1 \\
    0 & 0 & 1
\end{pmatrix} & q \text{ odd or } q\not\equiv -1\pmod{3},\\
    \begin{pmatrix}
    1 & \xi^3 & x_2 \\
    0 & 1 & \xi^{3q} \\
    0 & 0 & 1
\end{pmatrix} & q \text{ even and } q\equiv -1\pmod{3},
    \end{cases} \]
    with $x_2 \in \F_{q^2}^{\times}$.
\end{proposition}

\begin{proof}
It is sufficient to prove the first claim for $x \in S_1$ and $y \in S_2$. Hence, we can write 
\[
x = \begin{pmatrix}
    1 & x_1 & x_2 \\
    0 & 1 & x_3 \\
    0 & 0 & 1
\end{pmatrix}, \quad y = \begin{pmatrix}
    1 & 0 & 0 \\
    y_1 & 1 & 0 \\
    y_2 & y_3 & 1
\end{pmatrix}, \]
\[xy = \begin{pmatrix}
    1 + x_1 y_1 + x_2 y_2 & x_1 + x_2 y_3 & x_2 \\
    y_1 + x_3 y_2 & 1 + x_3 y_3 & x_3 \\
    y_2 & y_3 & 1
\end{pmatrix},
\]
for some $x_i, y_i \in \F_{q^2}^\times$, $1 \leq i \leq 3$.

Suppose $xy \in C_3^{(0,0)}$. Then $xy$ lies in a Sylow $p$-subgroup of $G$ different from $S_1$ and $S_2$. Hence, by the regularity of the action of $S_1$ on $\Syl_p(G) \setminus \{S_1\}$, $xy = a^{-1}ba$ for some $a \in S_1$ and $b \in S_2\setminus Z(S_2)$ of the form
\[
a = \begin{pmatrix}
    1 & a_1 & a_2 \\
    0 & 1 & a_3 \\
    0 & 0 & 1
\end{pmatrix} \neq 1, \qquad b = \begin{pmatrix}
    1 & 0 & 0 \\
    b_1 & 1 & 0 \\
    b_2 & b_3 & 1
\end{pmatrix},
\]
with $b_1 \neq 0 \neq b_3$.
Equivalently, we want to find a solution for the equation $axy = ba$ such that $x,y,b \in C_3^{(0,0)}$.
This equality holds if and only if the following equations are all satisfied:
\begin{equation*}
\begin{aligned}[c]
    1 + x_1y_1 + x_2y_2 + a_1y_1 + a_1x_3y_2 + a_2y_2 &= 1, \\
    x_1 + x_2y_3 + a_1 + a_1x_3y_3 + a_2 y_3 &= a_1, \\
    x_2 + a_1 x_3 + a_2 &= a_2, \\ 
    y_1 + x_3 y_2 + a_3 y_2 &= b_1, \\
    1 + x_3 y_3 + a_3 y_3 &= 1 + a_1 b_1. 
\end{aligned} \qquad
\begin{aligned}[c]
    x_3 + a_3 &= a_2 b_1 + a_3, \\
    y_2 &= b_2, \\
    y_3 &= a_1 b_2 + b_3, \\
    1 &= 1 + a_2 b_2 + a_3 b_3, \\
    \, & \, 
\end{aligned}
\end{equation*}

Recalling that $x_3 = -x_1^q$ and $x_2 + x_2^q + x_1^{q+1} = 0$, we find a solution to this system of equations in terms of $x_1, x_2$ and $b_1$:
\begin{equation}\label{eq:otherValues}
\begin{aligned}
   a_1 &= \frac{x_2}{x_1^q}, & a_2 &= - \frac{x_1^q}{b_1}, & a_3 &= - \frac{x_2^q}{x_1}, \\
         y_1 &= -b_1^q x_1^{q-1}, & y_2 &= \frac{b_1^{q+1} x_2^q}{x_1^{q+1}}, & y_3 &= b_1 x_1^{q^2-q}, \\
         \, & \, & b_2 &= \frac{b_1^{q+1} x_2^q}{x_1^{q+1}}, & b_3 &= -b_1^q.  
\end{aligned}
\end{equation}

Now we must find values of $x_1$, $x_2$ and $b_1$ such that they yield triangular unitary matrices $a,b,x,y$.
That is, for $\lambda \in \{a,b,x,y\}$, we must have $\lambda_3 = \,- \lambda_1^q$ and $\lambda_2~+~\lambda_2^q + \lambda_1^{q+1} = 0$. 
In addition, we need to check that $b,x,y\in C_3^{(0,0)}$.

First, we check the unitary conditions.
We will assume that $x$ is taken to be a unitary matrix.
Clearly, $\lambda_3 = - \lambda_1^q$ is always satisfied for $\lambda \in \{a,b,y\}$. So we focus on the equation $\lambda_2 + \lambda_2^q + \lambda_1^{q+1} = 0$.
\begin{enumerate}
    \item The matrix $a$ is unitary if and only if
    \begin{align} 
         - \frac{x_1^q}{b_1} + \Big( - \frac{x_1^q}{b_1} \Big)^q + \Big( \frac{x_2^q}{x_1}\Big)^{q+1} &= 0, \nonumber \\
        \Longleftrightarrow \frac{x_1^q}{b_1} + \frac{x_1}{b_1^q} &= \frac{x_2^{q+1}}{x_1^{q+1}}. \label{eq:unitary_condition}
    \end{align}
    \item The matrices $b$ and $y$ are unitary if and only if
    \begin{align*}
        \frac{b_1^{q+1} x_2^q}{x_1^{q+1}} + \Big( \frac{b_1^{q+1} x_2^q}{x_1^{q+1}} \Big)^q + b_1^{q+1} &= 0, \\
        \Longleftrightarrow \frac{b_1^{q+1} x_2^q}{x_1^{q+1}} + \frac{b_1^{q+1} x_2}{x_1^{q+1}} + b^{q+1} &= 0, \\
        \Longleftrightarrow b_1^{q+1} (x_2^q + x_2 + x_1^{q+1}) &= 0.
    \end{align*}
    This equality is satisfied if $x$ is unitary since then $x_2 + x_2^q + x_1^{q+1} = 0$.
\end{enumerate}

In conclusion, we have a solution for $axy=ba$, with $a,b,x,y$ unitary and $x,y,b \in~C_3^{(0,0)}$ if and only if we take the values as in Eq. (\ref{eq:otherValues}) and in addition the following equations hold:
\begin{equation}\label{eq:finalEquationsC300}
\begin{cases}
    x_2 + x_2^q + x_1^{q+1} = 0,\\
    \frac{x_1^q}{b_1} + \frac{x_1}{b_1^q} = \frac{x_2^{q+1}}{x_1^{q+1}},\\
    x,y,b \in C_3^{(0,0)}.
\end{cases}
\end{equation}
Therefore, we must find values of $x_1$, $x_2$, and $b_1$ such that Eq. (\ref{eq:finalEquationsC300}) holds.

\bigskip

\noindent
\textbf{Case 1.} Suppose $q$ is odd or $q \not\equiv -1 \pmod 3$. Then we can pick values of $x_1,x_2,b_1$ such that Eq. (\ref{eq:finalEquationsC300}) holds.

\begin{proof}[Proof of Case 1]

Set $x_1 = 1$. 
In order for $x$ to be unitary, $x_2$ then needs to satisfy $x_2 + x_2^q +1 = 0$. 
Furthermore, substituting $x_1=1$, Eq. (\ref{eq:unitary_condition}) becomes
\[
x_2^{q+1} = b_1^{-1} + b_1^{-q}.
\]
Dividing both sides by $-x_2^{q+1}$ and recalling that $(x_2^{q+1})^q = x_2^{q+1}$, this is equivalent to 
\[
(-b_1 x_2^{q+1})^{-1} + \big( (- b_1 x_2^{q+1})^{-1} \big)^q + 1 = 0.
\]
Hence, we can choose $(-b_1 x_2^{q+1})^{-1} = x_2$, implying $b_1 = -x_2^{-(q+2)}$.
Now we have that $x,y,a,b$ are unitary matrices, and by construction, they satisfy the equation $axy=ba$.
It remains to prove that $x,y,b\in C_3^{(0,0)}$.

First, if $q\not \equiv -1 \pmod{3}$, then any matrix
\[ z = \begin{pmatrix}
    1 & z_1 & z_2\\
    0 & 1 & z_3\\
    0 & 0 & 1
\end{pmatrix} \in S_1,\]
such that $z_1\neq 0$ automatically lies in $C_3^{(0,0)}$.
Similar conclusions hold if $z\in S_2$.
Thus, as $x_1 = 1$, $b_1\neq 0$ and $y_1 = -b_1^qx_1^{q-1}\neq 0$, we see that $x,y,b\in C_3^{(0,0)}$.
Hence these values of $x,y,a,b$ satisfy Eq. (\ref{eq:finalEquationsC300}) if $q\not\equiv -1\pmod{3}$.

Now, suppose that $q\equiv -1 \pmod{3}$, so per assumption $q$ must be odd.
By Lemma \ref{lem:PSU(3,q)_C3conjugates}, we need to show that $x_1,y_1,b_1$ are cubes.
Since $x_1=1$, $x_1$ is a cube.
Next, $y_1 = -b_1^qx_1^{q-1}$, and as $-1$ is a cube, we see that $y_1$ is a cube if $b_1$ is.
Note that $b_1 = -x_2^{-(q+2)}$, and $3\mid q+1$, so $3\nmid q+2$.
Hence $b_1$ is a cube if and only if $x_2$ is a cube, so it remains to check that we have a solution for $x_2+x_2^q+1=0$ where $x_2$ is a cube.
As $q$ is odd, here we can just take $x_2 = -\frac{1}{2}$.
This concludes the proof of this case.
\end{proof}

\noindent
\textbf{Case 2.} Suppose $q$ is even and $q \equiv -1 \pmod 3$. Then we can pick values of $x_1,x_2,b_1$ such that Eq. (\ref{eq:finalEquationsC300}) holds.

\begin{proof}[Proof of Case 2]
Suppose $b_1 = b_1^q \neq 0$, or equivalently, $b_1^{q-1} = 1$.
Then Eq. (\ref{eq:unitary_condition}) becomes
\[
\frac{x_2^{q+1}}{x_1^{q+1}} = \frac{1}{b_1} (x_1^q + x_1).
\]
This implies $b_1 = \frac{x_1^{q+1} (x_1 + x_1^q)}{x_2^{q+1}}$.
Then, $b_1^{q-1} = 1$ if and only if $x_1 + x_1^q \neq 0$.
As $x_1\neq 0$, the latter equation holds if and only if $x_1^{q-1} \neq 1$, that is, $x_1 \neq \xi^{m(q+1)}$ for any $0 \leq m < q-1$.
Recall that in addition we need $x_1, y_1$ and $b_1$ to be cubes and $x_2$ such that $x_2 + x_2^q + x_1^{q+1} = 0$.

Set $x_1 = \xi^3$. 
Then $x_1^{q-1}\neq 1$ since $q\neq 2$ by hypothesis, and clearly $x \in C_3^{(0,0)}$ as long as  $x_2 + x_2^q + x_1^{q+1} = 0$.
Similarly to Case 1, $y_1$ is a cube if and only if $b_1$ is a cube.
But, as $b_1^{q-1}=1$, we have that $b_1= \xi^{m(q+1)}$ for some $m\geq 0$, which is a cube since $3\mid q+1$.
All of this shows that, if $x_2 + x_2^q + x_1^{q+1} = 0$, then the given choices for $x_1,b_1$ satisfy Eq. (\ref{eq:finalEquationsC300}).
Therefore, it remains to find an appropriate value for $x_2$.
But we can always find such a value for $x_2$ as the trace map is surjective.
\end{proof}
These two cases complete the proof of Proposition \ref{prop:SU3_C3}.
\end{proof}

\begin{corollary} \label{cor:SU3(q)_qnot-1mod3}
Suppose $G= \SU_3(q)$, and $q \not\equiv -1 \pmod 3$, and let $\calC$ be the $G$-conjugacy class $C_3^{(0,0)}$.
Then $K_\calC$ is irreducible.
\end{corollary}

\begin{proof}
We follow Strategy 2 as outlined in Subsection \ref{sub:str2}, using a Sylow $p$-subgroup $S\leq G$ as $H$ there.
Let $\calG_S$ be the full subgraph of $\killingGraph{G}{\calC}$ with vertices $\calC \cap S$.

As $C_3^{(0,0)}$ intersects the centraliser of any element in $S$ by Lemma \ref{lem:PSU_3q_centC3}, in this case $K_\calC(x,y) \neq 0$ for any $x,y \in \calC \cap S$.
Hence, the cluster $\calG_S$ is connected.
Now, by Proposition \ref{prop:SU3_C3}, we see that there is another Sylow $p$-subgroup $S'\neq S$ such that there are edges connecting the vertices of the clusters $\calG_{S}$ with those of $\calG_{S'}$.
Finally, as $G$ is double transitive on the set of Sylow $p$-subgroups, we conclude by Strategy 2 that $\killingGraph{G}{\calC}$ is connected.
Therefore, $K_\calC$ is irreducible.
\end{proof}

Next, we want to prove that $K_{\calC}$ is irreducible for $\calC = C_3^{(0,l)}$, $0 \leq l \leq 2$, and $q\equiv -1\pmod{3}$.
However, in this case, it is not longer true that every centraliser of an element in a Sylow $p$-subgroup $S$ contains a representative of the class $C_3^{(0,l)}$, as now we have three conjugacy classes corresponding to the non-central elements of $S$, namely, $C_3^{(0,0)}, C_3^{(0,1)}$ and $C_3^{(0,2)}$. 
In fact, $K_{C_3^{(0,l)}}(x,y)$ might vanish for some elements $x,y \in C_3^{(0,l)} \cap S$ and our usual clusters $\calG_S$, with $S$ a Sylow $p$-subgroup, might not be connected.
To deal with this case, we will employ Strategy 3.

\begin{proposition} \label{prop:SU3(q)_q-1mod3}
    Let $G = \SU_3(q)$ for $q \equiv -1 \pmod 3$, $q \neq 2$. 
    Suppose $\calC~=~C_3^{(0,l)}$ for $0 \leq l \leq 2$.
    Then $K_\calC$ is irreducible.
\end{proposition}

\begin{proof}
We first look at $\calC = C_3^{(0,0)}$.

By Proposition \ref{prop:SU3_C3} and double transitivity, there are elements $y \in \calC \cap S_2$ and $x\in \calC\cap S_1$ such that $K_{\calC}(x,y)\neq 0$, and $x$ is as given in the Moreover part of that proposition.
Note that we can take $g \in \GU_3(q)$ to be a diagonal matrix such that $x^g \in S_1$ looks like one of the representatives of $C_3^{(0,0)}$ as given in Table \ref{tab:PSU(3,q)_conjclasses}.
Indeed, if $q$ is odd then $g=1$ and we set $u_0 := x$.
If $q$ is even, we can pick
\[ g = \begin{pmatrix}
    \xi^{3} & 0 & 0 \\
    0 & 1 & 0 \\
    0 & 0 & \xi^{-3q}
\end{pmatrix}, \]
so $u_0 := x^g$ looks like one of the representatives of $C_3^{(0,0)}$ as given in Table \ref{tab:PSU(3,q)_conjclasses}.

In any case, $g\in N_{\GU_{3}(q)}(S_2)$, so $Z(S_2)^g = Z(S_2)$.
Let $C_0$ be the connected component of the graph $\killingGraph{G}{\calC}$ that contains the vertices $x$ and $y$, and denote by $M$ its stabiliser under the $G$-action. 
Then we have $x\in M$ and $Z(S_2)\leq C_G(y) \leq M$.
Also, clearly we have $Z(\SU_3(q))\leq M$.
Conjugating by $g$ now gives
\[ M^g \geq \gen{x^g, Z(S_2)^g, Z(\SU_3(q))^g} = \gen{u_0, Z(S_2), Z(\SU_3(q))},\]
and the latter subgroup is $\SU_3(q)$ by Theorem \ref{thm:Nuzhin_gen_PSU3(q)}.
Therefore, $M = \SU_3(q)$ and $K_\calC$ is irreducible.

Now suppose $\calC = C_3^{(0,l)}$ with $l=1,2$. By Lemma \ref{lem:SU3(q)_C3_conjugates}, $\calC$ is a conjugate of $C_3^{(0,0)}$ in $\GU_3(q)$. In fact, all three conjugacy classes become a single class $\mathcal{D}$ in $\GU_3(q)$.
This implies that $\calG(\SU_3(q), \calC)$, as a subgraph of $\calG(\GU_3(q), \mathcal{D})$, is isomorphic to $\calG(\SU_3(q), C_3^{(0,0)})$, and thus also connected. 
As a result, $K_\calC$ is irreducible.
\end{proof}

Our results on conjugacy classes of unipotent elements in special unitary groups of dimension three are summarised in the following theorem. 

\begin{theorem} \label{thm:unitary_conclusion}
    Let $G$ be one of the groups $\SU_3(q)$ or $\PSU_3(q)$, with $q \neq 2$, and let $\calC$ be a conjugacy class of unipotent elements.
    Then $K_\calC$ is irreducible if and only if $\calC$ is not a conjugacy class of involutions.
    If $\calC$ is the class of involutions, then $K_{\calC}$ is non-degenerate.
\end{theorem}

\begin{proof}
    If $G = \SU_3(q)$, this follows from Proposition \ref{prop:PSU(3,q)_reduc}, Corollary \ref{coro:psu_C2_qodd_irred}, Corollary \ref{cor:SU3(q)_qnot-1mod3} and Proposition \ref{prop:SU3(q)_q-1mod3}. 
    
    For $G = \PSU_3(q)$ and $\calC$ a conjugacy class of non-involutive unipotent elements, this now follows by Lemma \ref{lm:connectedModCenter}.
    For $\calC$ the class of involutions when $q$ is even, this is Proposition \ref{prop:PSU(3,q)_reduc}.    
\end{proof}

\begin{remark}
    For the results in this section, we did not include the groups $\SU(3,2)$ and $\PSU(3,2)$.
    We verified computationally in GAP that in these groups, the conjugacy class of involutions also has a reducible and non-degenerate Killing form, and the three conjugacy classes of elements of order four have an irreducible Killing form.
\end{remark}


\section{Unipotent elements in Suzuki groups} \label{sec:suzuki}

In this section, we look at the unipotent elements in Suzuki groups $\Sz(q)$.
Recall that this group is always defined in characteristic two, so we are looking at $2$-elements. 
Our aim is to prove that $K_{\calC}$ is reducible and non-degenerate if $\calC$ is the conjugacy class of involutions, and irreducible for $\calC$ a conjugacy class of elements of order four.

The main results on conjugacy classes, Sylow 2-subgroups and the character table for the Suzuki groups were established in \cite{zbMATH03173723}. 
Below we recall the main facts needed here.

\begin{lemma} \label{lem:suzuki_facts}
    Let $G = \Sz(q)$ with $q = 2^{2n+1}$ for some $n \geq 1$. Let $r = 2^{n+1}$. 
    \begin{enumerate}
        \item $G$ has a unique conjugacy class of involutions of size $(q-1)(q^2+1)$, which we will denote as $X$.
        \item $G$ has two conjugacy classes of order-$4$ elements, each of size $\frac{(q^2-q)(q^2+1)}{2}$. Both classes are not real, and we will denote them by $Y$ and $Y^{-1}$ respectively.
        \item $G$ has $q^2 + 1$ Sylow $2$-subgroups, each of order $q^2$ and with normaliser of order $q^2(q-1)$.
        \item For each Sylow $2$-subgroup $S$, $Z(S)$ is elementary abelian of order $q$ and contains all the involutions of $S$. 
        \item $|S \cap X| = q-1$ and $|S \cap Y| = |S \cap Y^{-1}| = \frac{q(q-1)}{2}$.
        \item For $x \in S \cap Y$ or $x \in S \cap Y^{-1}$, $C_G(x)$ is completely contained in $S$ and has order $2q$.
        Moreover, $C_G(x) = \gen{x}Z(S)$.
        For $x\in Z(S)$ an involution, $C_G(x) = S$.
    \end{enumerate}
\end{lemma}
See Lemma 1 and Proposition 1 in \cite{zbMATH03173723} for more details.

In Table \ref{tab:Sz_characters} we recall part of the character table of $G$, as given in \cite{zbMATH03173723}.

\begin{table}[ht]
    \centering
    \begin{tabular}{c|c c c c c}
     & $\chi^1$ & $\chi^{q^2}$ & $\chi_k^{q^2+1}$ & $\chi_{k,s}^{(q-sr+1)(q-1)}$ & $\chi_k^{r(q-1)/2}$ \\
     & & & $k=1,\ldots, \frac{q}{2}-1$ & $k=1,\ldots, \frac{q+sr}{4}, s=\pm 1$ & $k=1,2$ \\
    \hline
    $1$ & $1$ & $q^2$ & $q^2+1$ & $(q-sr+1)(q-1)$ & $\frac{r(q-1)}{2}$ \\
    $X$ & $1$ & $0$ & $1$ & $sr-1$ & $- \frac{r}{2}$ \\
    $Y$ & $1$ & $0$ & $1$ & $-1$ & $-(-1)^k\frac{ri}{2}$ \\
    $Y^{-1}$ & $1$ & $0$ & $1$ & $-1$ &$(-1)^k \frac{ri}{2}$ \\
    \end{tabular}
    \caption{The irreducible characters of $\Sz(q)$ on conjugacy classes of unipotent elements.}
    \label{tab:Sz_characters}
\end{table}

We first look at the Killing form on the conjugacy class of involutions.

\begin{theorem} \label{thm:Sz_reduc}
    Let $G = \Sz(q)$ with $q = 2^{2n+1}$ and $n \geq 1$ and let $\calC = X$ be the conjugacy class of involutions. Then $\killingGraph{G}{\calC} = \commGraph{\calC}$, and $K_\calC$ is reducible and non-degenerate.
\end{theorem}

\begin{proof}
    Recall that $C_G(x) = S_x$ for any involution $x\in \calC$, where $S_x$ is the unique Sylow $2$-subgroup containing $x$.
    Hence, $C_G(xy)\cap \calC \neq \emptyset$, for $x,y\in \calC$, if and only if $xy$ is a $2$-element.
    We show that, for involutions $x,y\in \calC$, $xy$ is a $2$-element if and only if $S_x = S_y$, that is, they lie in the same Sylow $2$-subgroup.
    
    Suppose $xy$ has order four. Then $xyxy$ is an involution, contained in a Sylow 2-subgroup $S$, and $x,y \in C_G(xyxy) = S$. 
    As $x$ and $y$ are involutions, they must be central in $S$ and hence $xy\in Z(S)$ is an involution, a contradiction.
    
    On the other hand, $xy$ is an involution if and only if $x$ and $y$ commute and thus $x, y \in C_G(x) = S_x$ are contained in the same Sylow 2-subgroup.

    We conclude that $K_\calC(x,y) \neq 0$ if and only if $x$ and $y$ lie in the same Sylow $2$-subgroup.
    Hence $K_{\calC}$ is reducible, and $\killingGraph{G}{\calC} = \commGraph{\calC}$, with one connected component per Sylow $2$-subgroup.
    Moreover, each component is a complete graph on $q-1$ vertices.

    From the above computations, we can determine the values of $K_{\calC}(x,y)$ for $x,y \in \calC$:
    \begin{equation*}
    K_{\calC}(x,y) = \begin{cases}
        (q^2+1)(q-1) & x = y,\\
        q-1 & x,y \text{ lie in the same Sylow $2$-subgroup},\\
        0 & \text{otherwise.}
    \end{cases}
    \end{equation*}
    Hence, if we order the elements of $\calC$ by the Sylow $2$-subgroups, the Killing form $K_\calC$ becomes a block diagonal matrix,
    where each block on the diagonal has dimension $q-1$ and entries
    \[
    \begin{pmatrix}
        (q^2 + 1)(q-1) & q-1 & \dots & q-1 \\
        q-1 & (q^2 + 1)(q-1) & \dots & q-1 \\
        \vdots & \vdots & \ddots & \vdots \\
        q-1 & q-1 & \dots & (q^2 + 1)(q-1)
    \end{pmatrix}.
    \]

    As these matrices are invertible (by Miller's criterion with $E = q^2(q-1)I_{q-1}$ and $H = (q-1)\Theta_{q-1,q-1}$), we conclude that the Killing form $K_\calC$ is reducible and non-degenerate.
\end{proof}

Now we can conclude with the proof of Theorem \ref{thm:involutions_Lietype}.

\begin{proof}
[Proof of Theorem \ref{thm:involutions_Lietype}]
It follows from Example \ref{ex:PSL(2,q)_irred}, Propositions \ref{prop:psl(2,q)_reduc} and \ref{prop:PSU(3,q)_reduc}, and Theorem \ref{thm:Sz_reduc}.
\end{proof}

\smallskip

Next, we study the Killing form for the conjugacy classes of elements of order four.

\begin{theorem}
\label{thm:Sz_irred}
Let $G$ be the Suzuki group $\Sz(2^{2n+1})$, with $n \geq 1$, and let $\calC$ be a conjugacy class of elements of order four (that is, $\calC = Y$ or $Y^{-1}$). Then $K_\calC$ is irreducible.
\end{theorem}

\begin{proof}
We follow Strategy 2 from Subsection \ref{sub:str2} by taking $H$ a Sylow $2$-subgroup of $G$.
Recall that $G$ is double transitive on the set of its Sylow $2$-subgroups.

Let $x \in \calC$ and let $S$ be the unique Sylow $2$-subgroup that contains $x$.
As $C_G(x) = \gen{x}Z(S)$, this centralizer contains elements of the classes $X, Y$ and $Y^{-1}$.
Therefore, for $\calC \in \{Y, Y^{-1}\}$ and $x, y \in \calC$, $K_\calC(x,y) \neq 0$ if and only if $xy \in \{1\} \cup  X \cup Y \cup Y^{-1}$, that is, if and only if $xy$ is a $2$-element.
In particular, for any two $x,y \in \calC$ contained in the same Sylow $2$-subgroup, $K_\calC(x,y) \neq 0$.
This means that the full subgraph of $\killingGraph{G}{\calC}$ spanned by the elements of $\calC\cap S$, with $S$ a Sylow $2$-subgroup of $G$, is a complete graph (and in particular connected). 

Next, we use the character formula from Theorem \ref{thm:characterFormula} to count the number of pairs $(x,y) \in~\calC^2$ such that $xy$ is a $2$-element.
Note that $xy$ cannot be the neutral element as $\calC$ is not real.
Without loss of generality, we may assume $\calC = Y$.
By Table \ref{tab:Sz_characters}, we get the following:

\begin{align*}
    \Phi_G(Y,Y,X) &= |\{ (x, y) \in Y^2 : xy \in X \}| \\
    &= \frac{|Y|^2 |X|}{|G|} \sum_{\chi \in \Irr(G)} \frac{(\chi(Y))^2 \overline{\chi(X)}}{\chi(1)} \\
    &= \frac{q^2 (q^2+1)(q-1)(2q+1)}{4}.
\end{align*}
Similarly, we have
\begin{align*}
    \Phi_G(Y,Y,Y)  = \Phi_G(Y,Y,Y^{-1})
    &= \frac{q^2 (q^2+1)(q-1)(q^2-q-2)}{8}.
\end{align*}
Thus, in total, there are 
\[
\Phi_G(Y,Y,X) + \Phi_G(Y,Y,Y) + \Phi_G(Y,Y,Y^{-1}) = \frac{q^2(q^2+1)(q-1)(q^2+q-1)}{4}
\]
pairs of elements $(x,y) \in Y^2$ such that $K_Y(x,y) \neq 0$.

Next, we count how many of these pairs consist of two elements $x,y$ that are contained in the same Sylow $2$-subgroup $S$.
Since for any two elements $x,y\in Y\cap S$ we have $K_Y(x,y)\neq 0$,
we have
\[ |Y \cap S|^2 |G:N_G(S)| = \frac{q^2(q^2+1)(q-1)^2}{4} \]
pairs $(x,y)\in Y^2$ such that $x,y$ lie in the same Sylow $2$-subgroup.

Aside from these pairs, we get
\[
\frac{q^2(q^2+1)(q-1)(q^2+q-1)}{4} - \frac{q^2(q^2+1)(q-1)^2}{4} = \frac{q^4(q^2+1)(q-1)}{4}
\]
pairs $(x,y)\in Y^2$ such that $x,y$ lie in distinct Sylow $2$-subgroups.
Hence, there are two different Sylow $2$-subgroups $S_1,S_2$ and elements $x\in Y\cap S_1$, $y\in Y\cap S_2$ such that $K_Y(x,y)\neq 0$.
By Strategy 2, we conclude that $\killingGraph{G}{Y}$ is connected, that is, $K_Y$ is irreducible.
\end{proof}


\section{Unipotent elements in Ree groups} \label{sec:ree}

In this section, we look at the unipotent elements of the Ree groups $\Ree(q)$, where $q$ is an odd power of $3$.
Our aim is to prove that conjugacy classes of unipotent elements, that is, classes of $3$-elements, have irreducible Killing forms.

We refer to \cite{Ree} for the main results on conjugacy classes and the character table for the Ree groups.
We recall below the main results that we will need from \cite{Ree}.

\begin{lemma}
Let $G = \Ree(q)$, with $q=3^{2n+1}$ and $m = 3^n$.
Fix a Sylow $3$-subgroup $S$ of $G$.
Then $N_G(S) = S \rtimes \langle JR\rangle$ where $|J|=2$, $|R|=(q-1)/2$, and $J,R$ commute.
Moreover, we have a normal series $1\leq Z(S) \normal  P_1 \normal S$ such that:
\begin{enumerate}
    \item $Z(S)$ is elementary abelian of order $q$, $J$ acts by inversion on $Z(S)$, and $JR$ acts regularly on $Z(S)$. The elements of $Z(S)\setminus\{1\}$ belong to the same (real) conjugacy class, denoted by $X$.
    \item $P_1 = \Phi(S) = [S,S] = Z(S) \times F$ is elementary abelian of order $q^2$, with $R$ acting freely on $F$, and $J$ centralises $F$.
    In $P_1 \setminus Z(S)$, we have two conjugacy classes denoted by $T,T^{-1}$ (one is the inverse of the other, so they are not real classes).
    Note that $C_J(P_1) = F$.
    \item $S$ has order $q^3$ and hence $P_1$ has index $q$ in $S$. Moreover, every element $x\in S\setminus P_1$ has order nine and $x^3\in Z(S)$. In $S\setminus P_1$ we have exactly three conjugacy classes of the same size, which can be obtained as $Y$, $YT$ and $YT^{-1} = (YT)^{-1}$. The class $Y$ is real.
    \item In particular, $S/Z(S)$ is elementary abelian of order $q^2$.
    \item If $x\in S\setminus P_1$, then the elements of $x Z(S)$ are conjugates of $x$.
    \item $C_G(J) = \langle J\rangle\times \PSL_2(q)$.
\end{enumerate}
\end{lemma}

In Table \ref{tab:ree}, we summarise the information we need on the conjugacy classes of unipotent elements.
By a slight abuse of notation, we have denoted by $J$ the class of the element $J$ described above.

\begin{table}[ht]
    \centering
    \begin{tabular}{c|c|c|c|c|c}
       Class & $|x|$ & $|\calC|$ & $|\calC\cap S|$ & $C_G(x)$ & Classes in $C_G(x)$ \\
       \hline
       $X$  &  $3$ & $(q-1)(q^3+1)$ & $q-1$ & $S$ & $1,X,T^{\pm 1},Y, (YT)^{\pm 1}$\\
       $T$, $T^{-1}$ & $3$ & $\frac{q(q-1)(q^3+1)}{2}$ & $\frac{q(q-1)}{2}$ & $P_1\langle J\rangle$ &  $1,X,T^{\pm 1},J$\\
       $Y$ & $9$ & $\frac{q^2(q-1)(q^3+1)}{3}$ & $\frac{q^2(q-1)}{3}$ & $\langle x\rangle Z(S)$ & $1,X,Y$\\
       $YT$, $YT^{-1}$ & $9$ & $\frac{q^2(q-1)(q^3+1)}{3}$ & $\frac{q^2(q-1)}{3}$ & $\langle x\rangle Z(S)$ & $1,X,(YT)^{\pm 1}$
    \end{tabular}
    \caption{List of conjugacy classes $\calC$ of non-trivial unipotent elements in $\Ree(q)$. Here $x\in \calC\cap S$ is any element, with $x\in F$ in case $\calC  =T$ or $T^{-1}$.
    We always have $C_G(x) = C_{N_G(S)}(x)$.}
    \label{tab:ree}
\end{table}

Next, we prove:

\begin{theorem}
\label{thm:Ree_irred}
Let $G$ be the Ree group $\Ree(3^{2n+1})$, with $n\geq 1$, and let $\calC$ be a class of non-trivial unipotent elements in $G$.
Then $K_{\calC}$ is irreducible.
\end{theorem}

\begin{proof}
We follow Strategy 2 when possible, and Strategy 3 otherwise (see Subsections \ref{sub:str2} and \ref{sub:str3} respectively).
Recall that $G$ is double transitive on $\Syl_3(G)$.
Therefore, we aim to show that:
\begin{enumerate}[label=(\roman*)]
    \item Each cluster $\calG_S$ for $S$ a Sylow $3$-subgroup is connected.
    \item There are two distinct Sylow $3$-subgroups $S_1$, $S_2$, and elements $x \in S_1 \cap \calC$, $y \in S_2 \cap \calC$ such that $K_\calC(x,y) \neq 0$.
\end{enumerate}
To prove (ii), we count the pairs of elements $(x,y) \in \calC^2$ such that $K_\calC(x,y) \neq 0$ and then show that this number is strictly larger than the number of pairs of elements lying in the same Sylow $3$-subgroup.
For $x \in \calC$, let $\calC_i$, $1 \leq i \leq r$, be the conjugacy classes that have non-empty intersection with $C_G(x)$ (see the last column of Table \ref{tab:ree}).
Note that $|\Syl_3(G)| = q^3+1$.
Therefore, we show
\begin{enumerate}
    \item[(ii')] $\sum_{i = 1}^r \Phi_G(\calC, \calC, \calC_i) > |\calC \cap S|^2 (q^3+1)$,
\end{enumerate}
with $S \in \Syl_3(G)$.
Recall that $\Phi_G(\calC, \calC, \calC_i)$ can be computed using Theorem \ref{thm:characterFormula} as explained in Section \ref{sec:overall}.
    
Now we explore each case individually.

\bigskip

\noindent
\textbf{Case 1.} The theorem holds for $\calC = X$.
\begin{proof}
Here (i) is immediate, and indeed $\calG_S$ is a complete graph with $q-1$ vertices.
On the other hand, using Table \ref{tab:ree} and Theorem \ref{thm:characterFormula}, we have
\[ \sum_{\calC_i = 1,X,Y,T^{\pm1}, (YT)^{\pm1}} \Phi_G(\calC, \calC, \calC_i) = (q-1)(q^3+1)(q^3+q-1).\]
This is strictly larger than $|X \cap S|^2 (q^3+1) = (q-1)^2(q^3+1)$, so (ii') holds.
This finishes the proof of this case.
\end{proof}

\noindent
\textbf{Case 2.} The theorem holds for $\calC = T$ or $T^{-1}$.
\begin{proof}
We proceed as before.
Item (i) is immediate, and $\calG_S$ is also a complete graph with $\frac{q(q-1)}{2}$ vertices.
Using Table \ref{tab:ree} and Theorem \ref{thm:characterFormula}, we get
\[ 
\sum_{\calC_i = 1,X,T^{\pm1}, J} \Phi_G(\calC, \calC, \calC_i) = \frac{1}{4} q^2 (q^3+1)(q^2-1)(2q-1).
\]
This is strictly larger than $|\calC \cap S|^2 (q^3+1) = \frac{1}{4}q^2(q^3+1)(q-1)^2$, establishing (ii').
\end{proof}

For the remaining cases, we are unable to prove item (i), but we can still show that (ii') holds. We will conclude that the Killing form is irreducible by using a maximal-subgroup argument.

\bigskip
\noindent
\textbf{Case 3.} The theorem holds for $\calC = Y$, $YT$ or $YT^{-1}$.
\begin{proof}
Here we can only say a priori that $\calG_S$ contains $\frac{q^2(q-1)}{3}$ vertices, but it is not necessarily connected.
Using Table \ref{tab:ree} and Theorem \ref{thm:characterFormula}, we get
\begin{align*}
    \sum_{\calC_i = 1,X,Y} \Phi_G(Y,Y, \calC_i) &= \frac{q^3(q^3+1)(q-1)}{27} (q^4 + q^3 + 10q^2 + 9), \\
    \sum_{\calC_i = 1,X,(YT)^{\pm1}} \Phi_G(\calC, \calC, \calC_i) &= \frac{q^4(q^3+1)(q-1)}{27} (2q^3 - q^2 + 8q - 12), \; \calC = (YT)^{\pm1}.
\end{align*}
Using the bound for the number of solutions, we see this is strictly bigger than 
\[
|Y \cap S|^2 (q^3+1) = |(YT)^{\pm1} \cap S|^2 (q^3+1) = \frac{q^4(q-1)^2}{9}(q^3+1),
\] 
establishing (ii').

This proves that there are different Sylow $p$-subgroups $S_1,S_2$ and elements $x~\in~\calC\cap S_1$ and $y\in\calC\cap S_2$ which form an edge in $\killingGraph{G}{\calC}$.
Let $M$ denote the stabiliser of the connected component containing $x$ and $y$.
Hence $C_G(x),C_G(y)\leq M$.
Note that $Z(S_1),Z(S_2)\leq M$, and $C_G(x)$ does not contain an involution.
Therefore, $M$ cannot centralise an involution.
By the table of maximal subgroups of the Ree groups (see \cite{ReeMxl}), as $q\neq 3$, we conclude that $M$ must be contained in a proper parabolic subgroup $P$, or else $M = G$.
If $M\neq G$, then $M\leq P$ normalises a unique Sylow $3$-subgroup, which contradicts the fact that $M$ contains the centre of two different Sylow $3$-subgroups.
Therefore, $M = G$ and by Strategy 3 from Subsection \ref{sub:str3}, $\killingGraph{G}{\calC}$ is connected.
This proves that $K_{\calC}$ is irreducible.
\end{proof}
We have shown that for any class of non-trivial unipotent elements, $K_{\calC}$ is irreducible.
The proof of the theorem is then concluded.
\end{proof}

\section{Symmetric and alternating groups} \label{sec:symmetric_alternating}

The goal of this section is to provide information on the irreducibility of the Killing form for symmetric and alternating groups, when we take conjugacy classes of involutions or of elements of order $p$ in degrees $p$ or $2p$.
We summarise the results of this section in the following theorem.

\begin{theorem}
\label{thm:summarySymn_Altn}
Let $G = \Sym_n$ or $\Alt_n$, $\calC$ a $G$-conjugacy class, and $x\in \calC$ a representative.
Let $s$ denote the number of fixed points of $x$ by the natural representation.
Then the following hold:
\begin{enumerate}
    \item If $\calC$ is a conjugacy class of involutions, then $K_{\calC}$ is reducible if and only if one of the following holds:
    \begin{enumerate}
        \item $G = \Alt_n$, $n\equiv 1\pmod{4}$ and $s = 1$.
        \item $G = \Sym_n$ and $s = 1$.
        \item $G = \Sym_n$ and $(n,s) = (4,2)$.
    \end{enumerate}    
    In such cases, $K_{\calC}$ is non-degenerate.
    \item If $n = p$ or $2p$ and $x$ is a $p$-element, then $K_{\calC}$ is irreducible.
\end{enumerate}
\end{theorem}

The following result was first established in \cite[Theorem 1.1]{MR1994533}.

\begin{proposition}
\label{prop:involutionsSn}
Let $n\geq 2$, and let $\calC$ denote the conjugacy class of an involution of $\Sym_n$.
Write $s$ for the number of fixed points of a representative of $\calC$.
Then $\commGraph{\calC}$, the commuting graph on $\calC$, is disconnected if and only if $s=1$, or $(n,s)=(4,2)$.
\end{proposition}

We obtain an analogous answer for alternating groups.

\begin{corollary}
\label{coro:involutionsAltn}
Let $n\geq 2$, and let $\calC$ denote the conjugacy class of an involution of $\Alt_n$.
Write $s$ for the number of fixed points of a representative of $\calC$.
Then $\commGraph{\calC}$ is disconnected if and only if $n\equiv 1 \pmod{4}$ and $s=1$.
\end{corollary}

\begin{proof}
Note that the conjugacy class $\calC_x$ of $x$ in $\Alt_n$ is the same as the conjugacy class of $x$ in $\Sym_n$.
Thus, by Proposition \ref{prop:involutionsSn}, we see that the commuting graph on $\calC_x$ is disconnected if and only if $n$ is odd and $x$ fixes only one point, or $n=4$ and $x$ fixes two points.
Note that the latter condition cannot hold since $x$ lies in $\Alt_n$.
Since $x$ is the product of an even number of transpositions, if $x$ has only one fixed point then we must have $n\equiv 1\pmod{4}$.
\end{proof}

Now we show that in the case the commuting graph of a conjugacy class of involutions $\calC$ is disconnected, $K_{\calC}$ is reducible and invertible.

\begin{proposition}
\label{prop:Symn_involution_reducible}
Suppose $n \geq 3$ is odd and let $x\in \Sym_n$ be an involution that fixes only one point, and let $\calC_x$ denote its $\Sym_n$-conjugacy class.
Then, as a matrix, $K_{\calC_x}$ is the scalar matrix $|\calC_x|$.
Thus, $K_{\calC_x}$ is reducible and non-degenerate.
Moreover, $\killingGraph{\Sym_n}{\calC_x} = \commGraph{\calC_x}$.

The same conclusion holds if $x\in \Alt_n$.
\end{proposition}

\begin{proof}
The conclusions for $\Alt_n$ follow easily from those for $\Sym_n$ since the $\Alt_n$-conjugacy class of $x$ equals the $\Sym_n$-conjugacy class.

Recall that $C_{\Sym_n}(x) \cong C_2 ^r\rtimes \Sym_r$, where $n=2r+1$ and $r+1$ is the number of orbits of $x$.
In particular, every element of $C_{\Sym_n}(x)$ fixes at least the point fixed by $x$.
Also note that $C_{\Sym_n}(x)\cap \calC_x = \{x\}$.
In fact, if $y\in \calC_x$ fixes the same point as $x$, say, without loss of generality, $n$, then $x,y\in \Sym_{n-1}$.
But $n-1$ is even and $x,y$ are $\Sym_{n-1}$-conjugate, so by Proposition \ref{prop:involutionsSn}, the commuting graph on the $\Sym_{n-1}$-conjugacy class of $x$ and $y$ is connected, contradicting that $x$ and $y$ lie in different connected components of $\commGraph{\calC_x}$.
Therefore, $x,y$ cannot fix the same point.

Next, we prove that if $y\in \calC_x$ is such that $xy$ commutes with some $z\in \calC_x$, that is, $K_{\calC_x}(x,y)\neq 0$, then $x = y$.
To this end, we argue by contradiction using induction on $n$.

Suppose that $x\neq y$ lie in $\calC_x$ and that $xy$ commutes with some $z\in \calC_x$.
Observe that we must have $y\neq z$.
By the previous paragraph, $y,z$ cannot fix the same point.
Let $k_1$ denote the point fixed by $z$.
Since $xy\in C_{\Sym_n}(z)$, we have $xy(k_1) = k_1$, and $k_1\neq y(k_1)=:m_1$.
Thus $(k_1 \ m_1)$ is a cycle of both $x$ and $y$.
We note also that $k_2:=z(m_1) \neq m_1$ since $k_1$ is the only point fixed by $z$ and $m_1\neq k_1$.
Next, $xy$ permutes the $z$-orbits, and since $\{k_2,m_1\}$ is a $z$-orbit with $xy(m_1)=m_1$, we have $xy(k_2)=k_2$.
Let $x'=(k_1\ m_1)x$, $y'=(k_1 \ m_1)y$, and $z'=(k_2\ m_1)z$, and let $S$ denote the symmetric group on $[n]\setminus \{k_1,m_1\}$, so $S\cong \Sym_{n-2}$.
Then $x',y',z'\in S$ and they are $S$-conjugate, where $z'$ fixes only the point $k_1$.
Also note that $xy=x'y'$, which clearly commutes with $z'$.
Therefore, by induction, we conclude that $x'=y'$.
But this then shows that $x=(k_1\ m_1)x' = y$, a contradiction.

Note that this proof shows that indeed $K_{\calC_x}(x,y) = 0$ if $x\neq y$, and it is $|\calC_x|$ otherwise.
This concludes the proof.
\end{proof}

Next, we consider symmetric and alternating groups of degree $p$, where $p$ is a prime, and look at conjugacy classes of non-trivial $p$-elements.

\begin{lemma}
\label{lm:connectedSp}
Let $p$ be a prime number, and $G=\Sym_p$ or $\Alt_p$. If $\calC$ is a conjugacy class of elements of order $p$ in $G$,
then $K_{\calC}$ is irreducible.
\end{lemma}

\begin{proof}
We assume first that $G = \Sym_p$ and apply Strategy 3 to show that $K_\calC$ is irreducible.
We can assume that $p\geq 5$ (the cases $p=2,3$ are easy to verify).
Let $c = (1,2,\ldots,p)$, $g = (p, p-1, p-2)$, and $d = c^g$.
Then
\[d = (1,2,\ldots,p-3,p,p-2,p-1),\]
and
\[ cd = (1,3,\ldots,p-4, p,2,4,\ldots,p-3,p-1,p-2).\]
So $cd$ is a cycle of length $p$, meaning that $c,d$ form an edge in the graph $\killingGraph{\Sym_p}{\calC}$.

Now, a Sylow $p$-subgroup of $\Sym_p$ is a conjugate of $S := \gen{c}$, and its normaliser $N_{\Sym_p}(S)$ is a Frobenius group $C_p\rtimes C_{p-1}$.
Indeed, this is a maximal primitive subgroup (cf. the O'Nan-Scott theorem).
Since $S$ is abelian, the elements of $S\cap \calC$ form a clique and thus lie in the same connected component of $\killingGraph{\Sym_p}{\calC}$.
Denote such a component by $C$.
In particular, since $N_{\Sym_p}(S)$ stabilises such a clique, it must stabilise $C$.

To conclude, note that, by the previous paragraph, the vertex $d$ also lies in $C$, and therefore the normaliser of the unique Sylow $p$-subgroup containing $d$ must also normalise such a component.
We see that there are two distinct maximal subgroups of $\Sym_p$ stabilising the component $C$, and so their span must be the whole group $\Sym_p$.
By Strategy 3, $\killingGraph{\Sym_p}{\calC}$ is connected.

Now suppose that $G = \Alt_p$.
Then $G$ has exactly two conjugacy classes of elements of order $p$ by the splitting criterion.
Note that these classes are permuted by $\Sym_p$, and hence they give rise to isomorphic graphs $\killingGraph{G}{\calC}$.
Hence, it is enough to establish the connectivity of the graph for one of these conjugacy classes. 
Next, note that the elements $c,d,g$ defined above lie in $\Alt_p$, so the same computation shows that, if $\calC$ is the $\Alt_p$-conjugacy class of $c$, then $c$ and $d=c^g$ lie in the same connected component with $\Alt_p$-stabiliser $M$.
The same argument as above shows that $M$ contains the normaliser (in $\Alt_p$) of the two (different and unique) Sylow $p$-subgroups containing $c$ and $d$ respectively, which are indeed maximal subgroups of $\Alt_p$.
Therefore, $M = \Alt_p$, and $\killingGraph{\Alt_p}{\calC}$ is connected.
\end{proof}

Let $\calG_1,\calG_2$ be two simple graphs (no loops or multiple edges).
The graph $\calG_1 \times \calG_2$ is the simple graph whose vertices are the pairs $(v_1,v_2)$ where $v_i$ is a vertex of $\calG_i$, and there is an edge $\{ (v_1,v_2), (w_1,w_2)\}$ if $\{v_1,w_1\}, \{v_2,w_2\}$ are edges of $\calG_1,\calG_2$ respectively, or, for some $i$, $v_i = w_i$ and $\{v_{3-i},w_{3-i}\}$ is an edge of $\calG_{3-i}$.
With this notation, we get the following corollary.

\begin{corollary}
\label{coro:connectedSpSp}
Let $p$ be a prime, and let $G_1,G_2$ be two copies of $\Sym_p$.
Let $G = G_1\times G_2$, and let $\calC$ be the set of elements of the form $c_1c_2$, where each $c_i\in G_i$ has order $p$.
Then $\calC$ is a conjugacy class of $G$ and $K_\calC$ is irreducible.
\end{corollary}

\begin{proof}
It is not hard to verify that $\calC$ is a $G$-conjugacy class.
Moreover, $\killingGraph{G}{\calC} \cong \killingGraph{G_1}{\calC_1}\times \killingGraph{G_2}{\calC_2}$, where $\calC_i$ is the $G_i$-conjugacy class of elements of order $p$. 
By a well-known result of Weichsel \cite[Theorem 1']{Weichsel}, this product (as we defined here) is a connected graph if and only if $\killingGraph{G_i}{\calC_i}$ is connected for $i=1,2$. 
Since $\killingGraph{G_i}{\calC_i}$ is connected by Lemma \ref{lm:connectedSp}, we conclude that $\killingGraph{G}{\calC}$ is connected.
\end{proof}

Now we look at the case of symmetric and alternating groups on $2p$ points, with $p$ a prime number.

\begin{proposition}
\label{prop:Alt2p_Sym2p}
    Let $G = \Sym_{2p}$ or $\Alt_{2p}$, $p \geq 5$ a prime, and let $\calC$ be a conjugacy class of non-trivial $p$-elements. Then $K_{\calC}$ is irreducible.
\end{proposition}

\begin{proof}
Let $c_1 = (1,\ldots,p)$ and $c_2=(p+1,\ldots, 2p)$.
Then there are two conjugacy classes of non-trivial $p$-elements in $\Sym_{2p}$: the class represented by $c_1$ and the class represented by $c_1c_2$. We denote the first $\Sym_{2p}$-class by $\calC_1$, and the second one by $\calC_2$.
By the splitting criterion for conjugacy classes in the alternating groups, we see that these classes do not split in $\Alt_{2p}$ because $c_1$ and $c_1c_2$ contain two cycles of the same length (two of length one for $c_1$, and two of length $p$ for $c_1c_2$).
This implies that $\calC_1,\calC_2$ constitute the $\Alt_{2p}$-conjugacy classes of non-trivial $p$-elements in $\Alt_{2p}$, and that $\killingGraph{\Alt_{2p}}{\calC_i} \cong \killingGraph{\Sym_{2p}}{\calC_i}$ for $i=1,2$.
Thus, for simplicity, we will work with the symmetric group, and write $G:= \Sym_{2p}$.

Let $S_1 := \gen{c_1,c_2}$, which is a Sylow $p$-subgroup of $G$, and it is elementary abelian of order $p^2$.
Let $\calG_i := \killingGraph{G}{\calC_i}$.
Our aim is to prove that $\calG_i$ is connected using Strategy 3 from Section \ref{sec:overall}.

Since any two elements in $S_1\cap \calC_i$ commute, the elements of $S_1\cap \calC_i$ give rise to a clique in $\calG_i$.
We denote by $N_i$ the stabiliser of the connected component of $\calG_i$ that contains the elements of $S_1\cap \calC_i$.
Following Strategy 3, we need to prove that $N_i = G$, which then implies that $\calG_i$ is connected.

Let $G_1 = \Sym_{ \{1,\ldots,p\} }$, $G_2 = \Sym_{ \{p+1,\ldots, 2p\} }$, and $t = (1,p+1)\ldots(p,2p)$.
Note that $[G_1,G_2] = 1$ and $t$ interchanges $G_1$ with $G_2$.
Then $M := (G_1 G_2)\gen{t} \cong \Sym_p \wr C_2$ is a maximal (imprimitive) subgroup of $\Sym_{2p}$.

\bigskip

\noindent
\textbf{Case 1.} We have that $N_1 = G$, so $\calG_1$ is connected.

\begin{proof}
First, we prove that $M$ is contained in $N_1$.
Note that $G_i$ centralises $c_{3-i}$, and hence $G_1G_2\leq N_1$.
Since $t$ interchanges $c_1$ and $c_2$, $t$ must also stabilise the connected component of $\calG_1$ containing $c_1,c_2$.
That is, $M\leq N_1$.

Next, we find an element $x\in G\setminus M$ such that $x\in N_1$. This will allow us to conclude that $M = N_1$ because $M$ is a maximal subgroup of $G$.

Let $g$ be the cycle $(2p,p,p-2)$.
Then
\[c_1^g = (1,2,\ldots,p-3,2p,p-1,p-2).\]
From this, it is not hard to show that
\[c_1c_1^g = (1,3,\ldots,p-4,2p,p-1,p,2,4,\ldots,p-5,p-3)\]
is a cycle of length $p$.
In particular, its centraliser contains an element of $\calC_1$.
Thus $c_1^g$ is a vertex in the connected component of $c_1$, and so $c_1^g\in N_1$.
Moreover, $c_1^g \notin M$, so $x = c_1^g$ satisfies the desired requirements, concluding the proof of this case.
\end{proof}

\noindent
\textbf{Case 2.} We have $N_2 = G$, so $\calG_2$ is connected.

\begin{proof}
First, it is not hard to show that $M\cap \calC_2 = (G_1\times G_2)\cap \calC_2$ is the $(G_1\times G_2)$-conjugacy class of $G_1\times G_2$ corresponding to $c_1c_2$.
By Corollary \ref{coro:connectedSpSp}, we see that $\killingGraph{G_1\times G_2}{M\cap \calC_2}$ is connected.
In particular, this implies that $M$ stabilises the connected component of $c_1c_2$ in $\calG_2$.
That is, $M\leq N_2$.

Now we proceed as in the previous case and show that there is an element $x\in N_2\setminus M$, forcing $N_2 = G$ by the maximality of $M$.
Indeed, in this case we can take $g = (p,2p)\notin M$ and 
\[ (c_1c_2)^g = (1,2,\ldots,p-1,2p)(p+1,\ldots,2p-1,p).\]
Then
\begin{align*}
    (c_1c_2) \cdot (c_1c_2)^g = \ & (1,3,\ldots,p-2,2p,p+2,p+4,\ldots,p+(p-1))\\
    & \cdot (2,4,\ldots,p-1,p+1,p+3,\ldots,2p-2,p)
\end{align*}
is an element of $\calC_2$, implying $K_{\calC_2}(c_1c_2, (c_1c_2)^g) \neq 0$.
Therefore, $x := (c_2c_2)^g\in N_2\setminus M$, finishing the proof of this case.
\end{proof}
These two cases have shown that $\killingGraph{\Sym_{2p}}{\calC_i} \cong \killingGraph{\Alt_{2p}}{\calC_i}$ is connected for $i=1,2$, as desired.
\end{proof}

\begin{proof}[Proof of Theorem \ref{thm:summarySymn_Altn}]
Suppose $\calC$ is the conjugacy class of an involution.
Then either the commuting graph $\commGraph{\calC}$ is connected (which implies that $K_{\calC}$ is irreducible), or, by Proposition \ref{prop:involutionsSn}, $\commGraph{\calC}$ is disconnected and either $(n,s) = (4,2)$ or $s = 1$.
The case $(n,s) = (4,2)$ is straightforward to verify.
If $s = 1$ (with $n\equiv 1\pmod{4}$ if $G =\Alt_n$ by Corollary \ref{coro:involutionsAltn}), $K_{\calC}$ is reducible and invertible by Proposition \ref{prop:Symn_involution_reducible}.

The cases where $n=p$ or $2p$ and $\calC$ is a conjugacy class of non-trivial $p$-elements follow from Lemma \ref{lm:connectedSp} and  Proposition \ref{prop:Alt2p_Sym2p} respectively.
\end{proof}

\section{Groups with a strongly \texorpdfstring{$p$}{p}-embedded subgroup} \label{sec:strongly_p_embedded}

In this section, we analyse the irreducibility of the Killing form on conjugacy classes of $p$-elements in simple groups with a strongly $p$-embedded subgroup and $p$-rank at least two.
These groups are exactly those whose commuting graph on subgroups of order $p$ is disconnected, and groups of Lie type and Lie rank one in characteristic $p$ belong to this family.
We refer the reader to \cite{GLS3} for general properties of simple groups.

A strongly $p$-embedded subgroup of $G$ is a subgroup $M$ that contains a Sylow $p$-subgroup of $G$ and such that $M\cap M^g$ has order prime to $p$ for every $g\in G\setminus M$.
On the other hand, $G$ is termed an almost simple group if its socle is a non-abelian simple group.
Equivalently, the generalised Fitting subgroup of $G$, denoted by $F^*(G)$, is a simple group.

The following theorem is part of the classification of finite simple groups, and can be found for instance in \cite[(6.1), (6.2)]{Aschbacher2} and \cite[Theorem 7.6.1]{GLS3}.
We write $C_p(G)$ for the commuting graph on subgroups of order $p$ of $G$.

\begin{theorem}
\label{thm:stronglypembedded}
Let $G$ be a finite group and let $p$ be a prime.
Then $C_p(G)$ is disconnected if and only if $G$ contains a strongly $p$-embedded subgroup.

If $G$ is almost simple, has $p$-rank at least two and contains a strongly $p$-embedded subgroup, then one the following holds:
\begin{enumerate}
    \item $F^*(G)$ is a simple group of Lie type and Lie rank one in characteristic $p$, with $|G:F^*(G)|$ prime to $p$,
    \item $p\geq 5$ and $F^*(G) \cong \Alt_{2p}$,
    \item $p = 3$ and $G \cong \Aut(\PSL_2(8))$ or $M_{11}$,
    \item $p = 3$, $F^*(G) \cong \PSL_3(4)$ and $G/F^*(G)$ is a $2$-group,
    \item $p = 5$ and $G \cong \Aut(\Sz(32))$,
    \item $p = 5$ and $F^*(G) \cong {}^2F_4(2)'$, McL, or Fi$_{22}$,
    \item $p = 11$ and $G = F^*(G) \cong J_4$.
\end{enumerate}
\end{theorem}

Note that if $|G:F^*(G)|$ is prime to $p$, then $C_p(G) = C_p(F^*(G))$.

\begin{theorem}
\label{thm:stronglyPEmbeddedIrreducible}
Let $G$ be a simple group of $p$-rank at least two such that $C_p(G)$ is disconnected.
Let $\calC$ be a $G$-conjugacy class of non-trivial $p$-elements.
Then $K_{\calC}$ is reducible if and only if $p = 2$, $G$ is a group of Lie type and Lie rank one in characteristic two, and $\calC$ is the unique conjugacy class of involutions.
\end{theorem}

\begin{proof}
By Theorem \ref{thm:stronglypembedded}, $G$ is one of the groups in cases (1-7) there.
In case (1), the conclusions hold by the results of the previous sections, namely,
Proposition \ref{prop:psl(2,q)_reduc}, Corollary \ref{coro:linear_irred}, and Theorems \ref{thm:unitary_conclusion}, \ref{thm:Sz_reduc}, \ref{thm:Sz_irred} and \ref{thm:Ree_irred}.

If case (2) holds, then the result follows from Proposition \ref{prop:Alt2p_Sym2p}.

For the cases in items (3-7), we check using GAP that the corresponding graph $\killingGraph{G}{\calC}$ is connected.
This is easily done for the cases in items (3) and (4) and the Tits group in item (6).
Note that the cases $G = \Aut(\PSL_2(8))$ with $p=3$, and $G = \Aut(\Sz(32))$ with $p=5$ do not arise as $G$ is simple.
For the sporadic groups McL ($p=5$), $Fi_{22}$ ($p=5$), and $J_4$ ($p=11$), we apply Strategy 3 from Section \ref{sec:overall} aided with GAP computations.
In these cases, we first note that the elements of $\calC$ that lie in the same Sylow $p$-subgroup $S$ form a clique in the graph $\killingGraph{G}{\calC}$.
In particular, the normaliser of $S$ must be contained in the stabiliser of the connected component $C$ of the elements of $S \cap \calC$.
Therefore, as seen in Subsection \ref{sub:str3}, for an element $x\in S\cap \calC$, $N_G(\langle x \rangle)$ must also stabilise the connected component.
Next, we find an element $y$ in $\calC\setminus S$ in the same connected component $C$, that is, for which there is $x\in \calC\cap S$ such that $K_\calC(x,y) \neq 0$.
To conclude, we show that the normalisers of $\gen{x}$ and $\gen{y}$ span the whole group. This is easy to verify if, for instance, such normalisers are maximal subgroups (this holds for $p=11$ in $J_4$ for the class of central elements of order $11$, labelled 11A, see \cite[Table 5.3i]{GLS3}).
\end{proof}

\section{Non-degeneracy in dihedral groups}
\label{sec:dihedral}

In this final section, we show that dihedral groups of order $2n$, with $n$ odd, are strongly non-degenerate (see \cite[p.10]{MR3682632}).
Recall that a finite group $G$ is strongly non-degenerate if for every $G$-stable subset $\calC\subseteq G\setminus \{1\}$ that generates $G$ and consists of real elements, the associated Killing form $K_{\calC}$ is non-degenerate.

\begin{theorem}
\label{thm:dihedralgroup}
Let $G= D_{2n}$ be the dihedral group of order $2n$, with $n$ odd, and let $\calC \subseteq G \setminus \{1\}$ be a $G$-stable generating subset.
Then $K_{\calC}$ is non-degenerate.
\end{theorem}

\begin{proof}
Let $\calC \subseteq G \setminus \{1\}$ be a $G$-stable set that generates $G$.
Note that every conjugacy class is real, and that $\calC$ must contain the conjugacy class of involutions.

Write $G = \gen{s,t : t^n = 1 = s^2, \ sts = t^{-1}}$.
Then there exist $m\geq 0$ and $1\leq j_1 < \cdots < j_m < n/2$ such that
\[ \calC = \{ t^{j_i}, t^{-j_i} : 1\leq i\leq m\} \cup \{ t^i s :  0\leq i < n\}.\]
We compute first the values of $K_{\calC}$ on $\calC$.
Write $N = |\calC| = 2m + n$.
We have:
\begin{align*}
    K(t^{j_i}, t^{j_i}) & = |\gen{t}\cap \calC| = 2m,\\
    K(t^{j_i}, t^{-j_i}) & = |G\cap \calC| = N,\\
    K(t^{\pm j_i}, t^{\pm j_k}) & = |\gen{t}\cap \calC| = 2m & (i\neq k),\\
    K(t^{j_i}, t^ks) & = |C_G(t^{j_i+k}s)\cap \calC| = 1,\\
    K(t^is,t^js) & = |C_G(t^{i-j})\cap \calC| = 2m & (i\neq j),\\
    K(t^is,t^is) & = |G \cap \calC| = N.\\
\end{align*}
Now choose the ordered basis
\[ t^{j_1}, \ldots, t^{j_m}, t^{-j_m}, \ldots, t^{-j_1}, s, ts, \ldots, t^{n-1}s,\]
and note that $K_{\calC}$ takes the following format over this basis:
\[ K = \left( \begin{matrix}
    2m \Theta_{2m,2m} + n \overline{I}_{2m,2m} & \Theta_{2m,n}\\
    \Theta_{n,2m} & 2m\Theta_{n,n} + n I_{n,n}
\end{matrix}\right).\]
Recall that $\Theta_{a,b}$ is the matrix of size $a\times b$ where all its entries are $1$, $I_{a,b}$ is the matrix of size $a\times b$ with ones on its main diagonal, and $\overline{I}_{a,a}$ is the matrix of size $a\times a$ that has ones on its main anti-diagonal.

Next, we show that $K$ is invertible, and for that, we use Schur's complement formula together with Miller's criterion (see Remark \ref{miller}).
Let
\[ A = 2m \Theta_{2m,2m} + n \overline{I}_{2m,2m}, \ \ B= \Theta_{2m,n},\]
\[ C = \Theta_{n,2m}, \ \ D = 2m\Theta_{n,n} + n I_{n,n}.\]
By Miller's criterion with $E = nI_{n,n}$ and $H = 2m \Theta_{n,n}$, $D$ is invertible.
Moreover, its determinant is
\[ \det(D) = n^{n-1}(2m+n+2m(n-1)) = n^n(2m+1).\]
Thus, by Schur's complement formula, we have
\[ \det(K) = \det(D) \cdot \det(A - BD^{-1}C).\]
By Remark \ref{miller}, the inverse of $D$ is given by
\[ D^{-1} = E^{-1} - \frac{1}{1+g}E^{-1}HE^{-1},\]
where
\[ g = \Tr(HE^{-1}) = \frac{2m}{n} \Tr(\Theta_{n,n}) = 2m. \]
Therefore, we see that
\[ D^{-1} = \frac{1}{n} I_{n,n} - \frac{1}{2m+1}\frac{2m}{n^2}\Theta_{n,n}.\]
Now we compute:
\begin{align*}
    A - BD^{-1}C & = 2m\Theta_{2m,2m} + n \overline{I}_{2m,2m} - \Theta_{2m,n} \left(\frac{1}{n} I_{n,n} - \frac{1}{2m+1}\frac{2m}{n^2}\Theta_{n,n} \right) \Theta_{n,2m}\\
    & = 2m \Theta_{2m,2m}+n \overline{I}_{2m,2m} + \left(\frac{2m}{2m+1}-1\right) \Theta_{2m,2m}\\
    & =\left( 2m - \frac{1}{2m+1}\right) \Theta_{2m,2m} + n\overline{I}_{2m,2m}\\
    & = \frac{2m(2m+1)-1}{2m+1} \Theta_{2m,2m} + n\overline{I}_{2m,2m}.
\end{align*}
The above matrix, which we shall call $A'$, is invertible since its determinant is given by
\[ \det(A') = \frac{n^n}{a^n} (a-1)^{2m-1}(a+(2m-1)) (-1)^m,\]
where $a = \frac{(2m+1)n}{2m(2m+1)-1}$.
Hence $K$ is invertible and
\[ \det(K) = n^n(2m+1) \frac{n^n}{a^n} (a-1)^{2m-1}(a+(2m-1)) (-1)^m.\]
This concludes the proof of the non-degeneracy of $K_{\calC}$ for $G$-stable generating subsets $\calC\subseteq G \setminus \{1\}$, where $G=D_{2n}$ and $n$ is odd.
\end{proof}

\bibliographystyle{abbrv} 
\bibliography{refs.bib}

\end{document}